\documentclass{amsart}
\usepackage{amssymb}
\usepackage{amsfonts}
\usepackage{amstext}
\usepackage{algorithmic}
\usepackage{algorithm}
\usepackage{graphicx}
\usepackage{epstopdf}
\usepackage[all]{xy}

\numberwithin{equation}{section}
\newtheorem{theorem}{Theorem}[section]
\newtheorem{lemma}[theorem]{Lemma}
\newtheorem{proposition}[theorem]{Proposition}
\newtheorem{corollary}[theorem]{Corollary}

\newtheorem{fact}[theorem]{Fact}
\theoremstyle{definition}
\newtheorem{definition}[theorem]{Definition}
\newtheorem{example}[theorem]{Example}

\theoremstyle{remark}
\newtheorem{remark}[theorem]{\bf{Remark}}

\newtheorem{convention}[theorem]{\bf{Conventions}}
\newtheorem{notation}[theorem]{\bf{Notation}}

\newcommand{\Sym}{{\rm{Sym}}}
\newcommand{\Star}{{\rm{Star}}}

\newcommand{\Aut}{{\rm{Aut}}}

\newcommand{\Nor}{{\rm{Nor}}}

\newcommand{\Acal}{{\mathcal{A}}}

\newcommand{\Ecal}{{\mathcal{E}}}
\newcommand{\Gcal}{{\mathcal{G}}}

\newcommand{\Lcal}{{\mathcal{L}}}
\newcommand{\Rcal}{{\mathcal{R}}}

\newcommand{\Vcal}{{\mathcal{V}}}
\newcommand{\Ocal}{{\mathcal{O}}}

\newcommand{\la}{{\triangleright}}
\newcommand{\ra}{{\triangleleft}}
\begin{document}

\title[Set theoretic solutions of YBE and graphs]{Set theoretic solutions of the Yang-Baxter equation, graphs and computations}
\thanks{The first author was partially supported by Grant SAB2005-0094 of the Ministry of
Education and Culture of Spain, the Royal Society, UK, and by Grant
MI 1503/2005 of the Bulgarian National Science Fund of the Ministry
of Education and Science. Collaboration was during a visit by both authors to the Isaac Newton Institute, Cambridge, where the work originates.}

\author{Tatiana Gateva-Ivanova and Shahn Majid}
\address{TGI: Institute of Mathematics and Informatics\\
Bulgarian Academy of Sciences\\
Sofia 1113, Bulgaria\\
S.M:  Queen Mary, University of London\\
School of Mathematics, Mile End Rd, London E1 4NS, UK}

\email{tatianagateva@yahoo.com, tatyana@aubg.bg,
s.majid@qmul.ac.uk}

\date{\today}

\begin{abstract}
We extend our recent work on set-theoretic solutions of the Yang-Baxter or braid relations with
new results about their automorphism groups, strong twisted unions of solutions and multipermutation solutions. We introduce and study graphs of solutions  and use our graphical methods
for the computation of solutions of finite order and their
automorphisms. Results include a detailed study
of solutions of multipermutation level 2.
\end{abstract}
\maketitle

\section{Introduction}

 It is well-known that certain matrix solutions of the braid or Yang-Baxter equations lead to braided
 categories, knot invariants, quantum groups and other important constructions, see  \cite{Ma:prim} for an introduction. However, these equations are also very interesting at the level of set maps $r:X\times X\to X\times X$ where $X$ is a set and $r$ is a bijection, a line of study proposed in \cite{Dri,Wei}. Solutions extend linearly to very special linear solutions but also lead to a great deal of combinatorics and to algebras with very nice homological properties including those relating to  the existence of noncommutative Groebner bases. We start with the definitions, using conventions from recent works \cite{T04,TSh}.

\begin{definition}
Let $X$ be a nonempty set and let $r: X\times X \longrightarrow
X\times X$ be a bijective map. We use notation $(X,r)$ and refer
to it as a \emph{quadratic set}. The image of $(x,y)$ under $r$ is
presented as
\begin{equation}
\label{r} r(x,y)=({}^xy,x^{y}).
\end{equation}
The formula (\ref{r}) defines a ``left action" $\Lcal: X\times X
\longrightarrow X,$ and a ``right action" $\Rcal: X\times X
\longrightarrow X,$ on $X$ as:
\begin{equation}
\label{LcalRcal} \Lcal_x(y)={}^xy, \quad \Rcal_y(x)= x^{y},
\end{equation}
for all $ x, y \in X.$ The map  $r$ is \emph{nondegenerate}, if the
maps $\Rcal_x$ and $\Lcal_x$ are bijective for each $x\in X$.   In
this paper we shall always be interested in the case where $r$ is nondegenerate,  as
will be indicated. As a notational tool, we  shall often identify
the sets $X\times X$ and $X^2,$ the set of all monomials of length
two in the free semigroup $\langle X\rangle.$
\end{definition}
\begin{definition}
\label{defvariousr}
\begin{enumerate}
 \item
$r$ is \emph{square-free} if $r(x,x)=(x,x)$ for all $x\in X.$
\item
\label{YBE} $r$ is \emph{a set-theoretic solution of the
Yang-Baxter equation} or, shortly \emph{a solution}
 (YBE) if  the braid relation
\[r^{12}r^{23}r^{12} = r^{23}r^{12}r^{23}\]
holds in $X\times X\times X,$ where the two bijective maps
$r^{ii+1}: X^3 \longrightarrow X^3$, $ 1 \leq i \leq 2$ are
defined as  $r^{12} = r\times Id_X$, and $r^{23}=Id_X\times r$. In
this case some authors also call $(X,r)$   \emph{a braided set}.
If in addition $r$ is involutive $(X,r)$ is called \emph{a
symmetric set}.
\end{enumerate}
\end{definition}
To each quadratic set $(X,r)$  we
associate canonical algebraic objects  generated by $X$ and with
quadratic defining relations $\Re =\Re(r)$ defined by
\begin{equation}
\label{defrelations} xy=zt \in \Re(r),\quad
  \text{whenever}\quad r(x,y) = (z,t).
\end{equation}

\begin{definition}
\label{associatedobjects}  Let $(X,r)$ be a quadratic set.

(i) The semigroup
\[
S =S(X, r) = \langle X; \Re(r) \rangle,
\]
 with a
set of generators $X$ and a set of defining relations $ \Re(r),$
is called \emph{the semigroup associated with $(X, r)$}.

(ii) The \emph{group $G=G(X, r)$ associated with} $(X, r)$ is
defined as
\[
G=G(X, r)={}_{gr} \langle X; \Re (r) \rangle.
\]

(iii) For arbitrary fixed field $k$, \emph{the $k$-algebra
associated with} $(X ,r)$ is defined as
\begin{equation}
\label{Adef} \Acal = \Acal(k,X,r) = k\langle X \rangle/(\Re(r)).
\end{equation}

If $(X,r)$ is a solution, then $S(X,r)$, resp.\ $G(X,r)$, resp.\
$\Acal(k,X,r)$ is called the {\em{Yang-Baxter semigroup}}, resp.\
the {\em{Yang-Baxter group}}, resp.\ the {\em{Yang-Baxter
algebra}} associated to $(X,r)$.
\end{definition}

\begin{example}
\label{trivialsolex} For arbitrary nonempty set $X$, \emph{the
trivial solution} $(X, r)$ is defined as $r(x,y)=(y,x),$ for all
$x,y \in X.$ It is clear that $(X,r)$ is the trivial solution
\emph{iff} ${}^xy =y$, and $x^{y} = x,$ for all $x,y \in X,$ or
equivalently $\Lcal_x= id_X =\Rcal_x$ for all $x\in X.$ In this case
$S(X,r)$ is the free abelian monoid, $G(X,r)$ is the free abelian
group,  $\Acal(k,X,r)$ the algebra of commutative polynomials in
$X$.
\end{example}

The algebras $\Acal(X,r)$ provided new classes of Noetherian rings
\cite{T94,T96}, Gorentstein (Artin-Schelter regular) rings
\cite{T96Preprint,T00,T04s} and so forth. Artin-Schelter regular
rings were introduced in \cite{AS} and are of particular interest.
The algebras $\Acal(X,r)$ are similar in spirit to the quadratic
algebras associated to linear solutions, particularly studied in
\cite{Man}, but have their own remarkable properties in the set
theoretic case. The semigroups $S(X,r)$ were studied particularly
in \cite{TSh} with a systematic theory of `exponentiation' from
the set to the semigroup by means of the `actions'
$\Lcal_x,\Rcal_x$ (which in the process become a matched pair of
semigroup actions) somewhat in analogy with the Lie theoretic
exponentiation in \cite{Ma:mat}.  Also in \cite{TSh} are first
results about construction of extensions of solutions that we
shall take further. There are many other works on set-theoretic
solutions and related structures, of which a relevant selection
for the interested reader is \cite{ESS,TM,GJO,JO,LYZ,Rump, Tak}.
Of particular interest in the present paper, we define:

\begin{definition}
When $(X,r)$ is nondegenerate,   $\Lcal:G(X,r) \longrightarrow
\Sym(X)$ is a group homomorphism defined via the left
action. We denote by $\Gcal= \Gcal(X,r)$ the subgroup
$\Lcal(G(X,r))$ of $\Sym(X)$.
\end{definition}
For example, $\Gcal(X,r) = \{id_X \}$ in the case of the trivial solution.

\begin{remark} \label{ybe} As is well known, see for example \cite{TSh},
a quadratic set   $(X,r)$ is a braided set (i.e. $r$ obeys the YBE)
 {\em iff} the following conditions hold
\[
\begin{array}{lclc}
 {\bf l1:}\quad& {}^x{({}^yz)}={}^{{}^xy}{({}^{x^y}{z})},
 \quad\quad\quad
 & {\bf r1:}\quad&
{(x^y)}^z=(x^{{}^yz})^{y^z},
\end{array}\]
 \[ {\rm\bf lr3:} \quad
{({}^xy)}^{({}^{x^y}{(z)})} \ = \ {}^{(x^{{}^yz})}{(y^z)},\]
 for all $x,y,z \in X$.
Clearly, conditions {\bf l1}, respectively, {\bf r1} imply that
the group $G(X,r)$ acts on the left, respectively on the right, on
the set $X$. This is needed for the definition of $\Gcal(X,r)$ above to make sense. \end{remark}

We also have some technical conditions and special cases of interest, particularly  certain `cyclicity' conditions. Systematic treatments were provided in \cite{T04,TSh} and we recall some essentials needed for the paper in the remainder of this introduction.

\begin{definition}
\label{extendedleftaction} \cite{TSh}
Given $(X,r)$ we extend the  actions ${}^x\bullet$ and  $\bullet
^x$ on $X$ to left and right actions on $X\times X$ as follows. For
$x,y,z \in X$ we define:
\[ {}^x{(y,z)}:=({}^xy,{}^{x^y}z), \quad
\text{and} \quad (x,y)^z:= (x^{{}^yz}, y^z).\] We say that $r$ is
respectively left and right invariant if
\[{\bf l2}: \; \; \;  \; \; \;
r({}^x{(y,z)})={}^x{(r(y,z))},\quad {\bf r2}: \; \; \; \; \; \;
r((x,y)^z)={(r(x,y))}^z\] hold for all $x,y,z\in Z$.
\end{definition}

\begin{definition}
\label{cyclicconditionsall} \cite{TSh},
A quadratic set $(X,r)$ is
called {\em cyclic} if the following conditions are satisfied
\[\begin{array}{lclc}
 {\rm\bf cl1:}\quad&  {}^{y^x}x= {}^yx \quad\text{for all}\; x,y \in
 X;
 \quad&{\rm\bf cr1:}\quad &x^{{}^xy}= x^y, \quad\text{for all}\; x,y \in
X;\\
 {\rm\bf cl2:}\quad
  &{}^{{}^xy}x= {}^yx,
\quad\text{for all}\; x,y \in X; \quad & {\rm\bf cr2:}\quad
&x^{y^x}= x^y \quad\text{for all}\; x,y \in X.
\end{array}\]
 We refer to these
conditions as {\em cyclic conditions}.
\end{definition}

 Theorem \ref{basictheorem} given below shows that every  square-free solution
 $(X,r)$ is cyclic, furthermore it satisfies condition \textbf{lri} and is a
 \emph{left and a right cycle set}, see the definitions below.

 \begin{definition}
\label{lri} \cite{TSh} Let $(X,r)$ be a quadratic set. We define the condition
condition
\[ \textbf{lri:}
\;\;\;\;\;\;\;\;\;\; ({}^xy)^x= y={}^x{(y^x)} \;\text{for all} \;
x,y \in X.\] In other words \textbf{lri} holds if and only if $(X,r)$ is nondegenerate and
$\Rcal_x=\Lcal_x^{-1}$ and $\Lcal_x = \Rcal_x^{-1}$
\end{definition}

\begin{remark}
\label{lrinondegenerate} By contrast with this involutiveness does not entail nondegeneracy. For
example \emph{the identity solution } $(X, r_X)$, with $r_{X}= id_{X\times X}$ is obviously involutive, but for each $X$ of
order $\geq 2$ it is not nondegenerate.
\end{remark}

 \begin{definition}
\label{csl&csr} \cite{Rump,TSh} Let  $(X,r)$ be
nondegenerate. $(X,r)$ is called \emph{a left cycle set} if
\[ \textbf{csl:} \; \; \;\; \; \;
{}{}^{({}^yt)}{({}^yx)}={}{}^{({}^ty)}{({}^tx)} \; \; \text{for
all}\; x,y,t \in X.\] The definition of  \emph{a right cycle set} is
analogous, see \cite{TSh}.
\end{definition}

\begin{proposition}
\label{lri&INVOLUTIVE} \cite[Proposition 2.24]{TSh} Let $(X,r)$ be
a quadratic set. Then any two of the following conditions imply the
remaining third condition.
\begin{enumerate}
\item
$(X,r)$ is involutive.
\item
$(X,r)$ is nondegenerate and cyclic, see Definition
\ref{cyclicconditionsall}.
\item
{\bf lri} holds.
\end{enumerate}
\end{proposition}
\begin{theorem}
\label{basictheorem} \cite[Theorem 2.34]{TSh} Suppose $(X,r)$ is
nondegenerate, involutive and square-free quadratic set. Then the
following conditions are equivalent:
\begin{enumerate}
\item \label{YBESquarefree} $(X,r)$ is a set-theoretic
solution of the Yang-Baxter equation. \item \label{YBEl1} $(X,r)$
satisfies {\bf l1}. \item \label{YBEl2} $(X,r)$ satisfies {\bf
l2}.
\item \label{YBEr1} $(X,r)$ satisfies {\bf r1}.
\item \label{YBEr2} $(X,r)$ satisfies {\bf r2}.
 \item\label{YBElr3}
$(X,r)$ satisfies {\bf lr3} \item \label{cycleset} $(X,r)$
satisfies {\bf csl}.
\end{enumerate}
In this case $(X,r)$ is cyclic and satisfies  {\bf lri}.
\end{theorem}
\begin{corollary}\label{constructivecor}
Every nondegenerate involutive square-free solution $(X,r)$ is uniquely determined by the
left action $\Lcal: X\times X \longrightarrow X,$ more precisely,
\[
r(x,y) = (\Lcal _x(y), \Lcal^{-1}_y(x)).
\]
Furthermore it  is cyclic.
\end{corollary}

\section{Homomorphisms of solutions, automorphisms}

\begin{definition}
\label{defhom} Let $(X,r_X)$ and $(Y, r_Y)$ be arbitrary solutions
(braided sets).
A map $\varphi: X \longrightarrow Y$ is a {\em{homomorphism of
solutions}}, if it satisfies the equality
\[(\varphi \times \varphi) \circ r_X =
r_Y\circ (\varphi \times \varphi).
\]
A bijective homomorphism of solutions  is called (as usual)
\emph{an isomorphism} of solution.
An isomorphism of the solution
$(X,r)$ onto itself is \emph{an r-automorphism}.

We denote by $Hom((X,r_X),(Y,r_Y))$ the set of all homomorphisms
of solutions $\varphi: X \longrightarrow Y$. The group of
$r$-automorphisms of $(X, r)$ will be denoted by $\Aut(X, r)$.
Clearly, $Aut(X,r)$ is a subgroup of $Sym(X).$
\end{definition}

\begin{remark} Let $(X, r_X), (Y, r_Y)$ be finite symmetric sets with \textbf{lri}.
Every homomorphism of solutions  $\varphi: (X, r_X)
\longrightarrow (Y, r_Y)$ induces canonically a homomorphism of
their graphs (see Definition \ref{defgraph}):
\[\varphi_{\Gamma}:
\Gamma(X, r_X) \longrightarrow \Gamma(Y, r_Y).\] Furthermore there
is a one-to one correspondence between $Aut(X,r)$ and
$Aut(\Gamma(X,r)),$ the group of \emph{automorphisms} of the
multigraph $\Gamma(X,r).$
\end{remark}

The next lemma is straightforward from the definition.
\begin{lemma}\label{leftright}
Let $(X,r_X)$ and $(Y,r_Y)$ be
solutions.
\begin{enumerate}
\item
 A map $\varphi: X \longrightarrow Y$ is a homomorphism of
solutions \emph{iff}:
\[
\varphi \circ \Lcal_x =\Lcal_{\varphi(x)} \circ\varphi {\rm{\ \ and\ \ }}
\varphi \circ \Rcal_x =\Rcal_{\varphi(x)} \circ \varphi \;\text{for all}\; x \in X.
\]
\item
If both $(X,r_X)$ and $(Y,r_Y)$ satisfy {\bf lri}, then $\varphi$ is
a homomorphism of solutions \emph{iff}
\[
\varphi \circ \Lcal_x =\Lcal_{\varphi(x)}\circ \varphi,  \;\text{for
all}\; x \in X.
\]
\item If $(X,r)$ obeys {\bf lri}, then
 $\tau \in Sym(X)$ is an automorphism of $(X,r_X)$ \emph{iff}
\begin{equation}
\label{autre}
 \tau\circ \Lcal_x  \circ \tau^{-1}=\Lcal_{ \tau(x)},
\;\text{for all}\; x \in X.
\end{equation}
\end{enumerate}
\end{lemma}

For example, if $(X,r)$ is the trivial solution then, clearly, $\Aut(X, r)=Sym(X).$ This is because
each $\Lcal_x=id_X$. More generally, (\ref{autre}) clearly implies:
\begin{corollary} The group $\Aut(X, r)$ is a subgroup of $\Nor_{\Sym(X)} \Gcal(X,r)$, the
normalizer of  $\Gcal(X,r)$ in $\Sym(X)$.
\end{corollary}

We now turn to two basic examples which will be used throughout
the paper. The first is a solution $(X,r)$ of order $6$, with
$\Gcal(X,r)\subsetneq Aut(X,r)$ for which $\Aut(X, r) \subsetneq
\Nor_{\Sym(X)}(\Gcal).$ We shall show in general in Section 4 that
$\Gcal(X,r)$ is a subgroup of $\Aut(X, r)$ \emph{iff} $(X,r)$ is a
multipermutation solution of level $2$, see
Theorem~\ref{mpl2theorem}.

\begin{figure}
\begin{center}
\includegraphics[scale=.8]{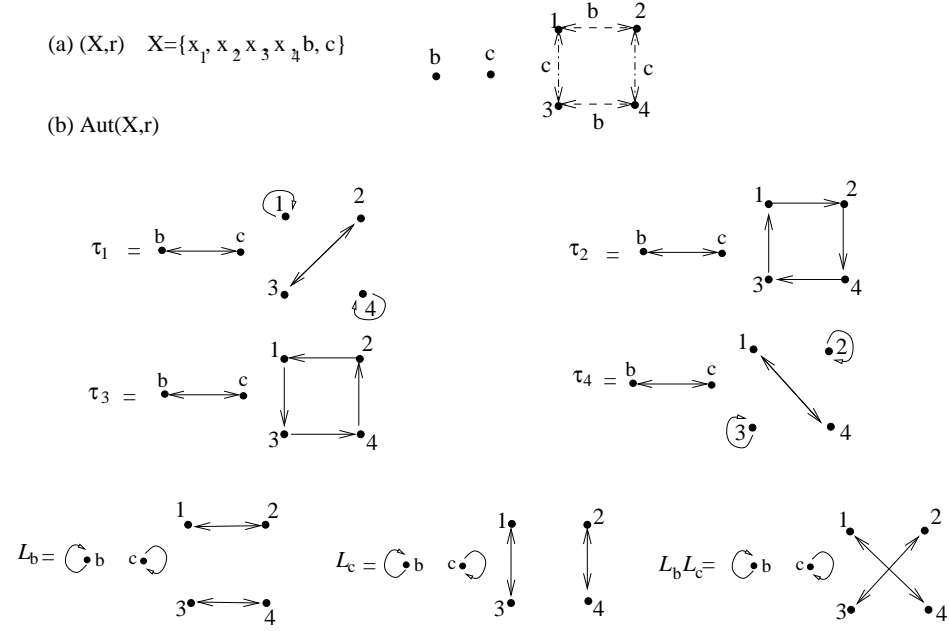}
\caption{(a) Graph of the solution of Example~\ref{autex1} and (b) its automorphism group.}
\label{figautex1}
\end{center}
\end{figure}

\begin{example}
\label{autex1} Let $(X,r)$ be the nondegenerate involutive square-free solution with
\[
X = \{x_1, x_2, x_3, x_4, b,c \},\quad \Lcal_b=(x_1x_2)(x_3x_4),\quad
\Lcal_c=(x_1x_3)(x_2x_4)
\]
and the remaining $\Lcal_{x_i}=id_X$. Here and elsewhere we use Corolary~\ref{constructivecor}.  Then $ \Gcal(X,r)=\langle \Lcal_b \rangle \times\langle \Lcal_c \rangle,$
so it is isomorphic to the {\em{Klein's group}} $\mathbb{Z}_2 \times
\mathbb{Z}_2$. Direct computations show that the set of
automorphisms consists of the following eight elements:
\[\begin{array}{lclc}
 id_X,  &\tau_1 = (bc)(x_2x_3), \quad &\tau_2= (bc)(x_1x_2x_4x_3),\\
\tau_3= (bc)(x_1x_3x_4x_2), \quad & \tau_4= (bc)(x_1x_4) , \quad
&\Lcal_b, \quad \Lcal_c, \quad\Lcal_b\circ\Lcal_c.\end{array}
\]
Furthermore, one has
\[
\tau_1^2 = 1,\quad\tau_2^4 =1,\quad \tau_1 \tau_2\tau_1^{-1}
=\tau_2^3 = \tau_3,\quad\tau_1\circ\tau_2 = \Lcal_c,\quad
\tau_2\circ\tau_1 =\Lcal_b,\quad\tau_1\circ\tau_2^2= \tau_4.
\]
It is easy to see that
\[
 Aut(X,r)={}_{gr} \langle \tau_1, \tau_2 \mid \quad \tau_1^2= 1,\tau_2^4 =1,
 \tau_1 \tau_2\tau_1^{-1} =\tau_2^3  \rangle \approx D_4.
\]
Next we show that $Aut(X,r)\subsetneq \Nor_{\Sym(X)}(\Gcal).$
Consider $\sigma=(x_1x_3)(x_2x_4)(bc)\in \Sym(X)$. One has:
\begin{equation}\label{eq01NorAut}
\sigma \circ \Lcal _{b} \circ \sigma ^{-1}=\Lcal _{b},\quad \sigma
\circ \Lcal _{c} \circ \sigma ^{-1}=\Lcal _{c},
\end{equation}
so $\sigma\in\Nor_{\Sym(X)}\Gcal.$ On the other hand $\sigma$ does
not satisfy the necessary condition (\ref{autre}) for being an
automorphism, so $\Aut(X, r)$ is a proper subgroup of
$\Nor_{\Sym(X)}(\Gcal).$ The  graph $\Gamma(X, r)$ and its automorphism group is given on figure \ref{figautex1}.
\end{example}

The second example will contain the above solution $(X,r)$ an $r$-invariant subset. (In fact $(X,r)$ itself contains  $X_0=\{x_1,x_2,x_3,x_4\}$ as an $r$-invariant subset.)

\begin{definition}
Let $(Z,r)$ be a quadratic set, $\emptyset \neq X\subseteq Z.$ $X$
is \emph{$r$-invariant} if $r(X\times X) \subseteq X\times X.$ In
this case we consider the restriction $r_X = r_{\mid X\times X}.$
Clearly, if  $(Z,r)$ obeys YBE then  $(X,r_X)$ is also a solution,
and inherits all ``good" properties of $(Z,r)$ like being non
degenerate, involutive, square-free, satisfying \textbf{lri }, the
cyclic conditions, etc.
\end{definition}

\begin{lemma}
\label{criterionLa} Let $(Z,r)$ be a nondegenerate
solution, with \textbf{lri} and $\alpha \in Z.$ Suppose $X$ is an
$r$-invariant subset of $Z$ and let $r_X$ be the restriction of $r$ on
$X\times X.$ Then

1) $\Lcal_{\alpha\mid X} \in \Aut(X, r_X)$ \emph{iff}
\begin{equation*}\label{LainAut}
{}^{\alpha^y}x={}^{\alpha}x  \quad \text{for all}\quad
x,y, \in X.
\end{equation*}
  In particular,
$\Lcal_{\alpha} \in \Aut(Z, r)$ \emph{iff} the displayed condition holds for
all $x,y \in Z.$

2) Furthermore, suppose that $r$ is involutive, then
$\Lcal_{\alpha\mid X} \in \Aut(X, r_X)$ \emph{iff}
\begin{equation*}
{}^{{}^y{\alpha}}x= {}^{\alpha}x \quad \text{for
all}\quad x,y, \in
X.
\end{equation*}
\end{lemma}
\begin{proof}
Let $\alpha\in Z$. By (\ref{autre}) $\Lcal_{\alpha\mid
X}$ is an
automorphism of $(X,r)$ \emph{iff}
\[
 \Lcal_{\alpha\mid X}\circ \Lcal_y
=\Lcal_{{}^{\alpha}y} \circ \Lcal_{\alpha\mid X},
\;\text{for all}\; y \in X
\]
or equivalently
\begin{equation}
\label{autre1} {}^{\alpha}{({}^yx)} = {}^{{}^{\alpha}
y}
{({}^{\alpha} x)}, \;\text{for all}\; x, y \in X
\end{equation}
By condition \textbf{l1} on $(Z,r)$ one has
\[
{}^{\alpha}{({}^yx)} = {}^{{}^{\alpha} y} {({}^{\alpha
^y} x)},
\;\text{for all}\; x, y \in X,
\]
which together with (\ref{autre1}) implies
\begin{equation}
\label{autre2} {}^{{}^{\alpha} y} {({}^{\alpha} x)}
 = {}^{{}^{\alpha} y} {({}^{\alpha ^y} x)},
\;\text{for all}\; x, y \in X
\end{equation}
By the non degeneracy of $(Z,r)$  (\ref{autre2}) holds
\emph{iff}
${}^{\alpha ^y} x = {}^{\alpha} x.$ This proves part
(1).

Assume now that  $r$ is involutive. Then $(Z,r)$ is a
symmetric set
so its retraction $([Z], r_{[Z]})$ is well defined,
see Definition~\ref{[x]def}.  By Lemma \ref{retractactionslemma1}
$([Z], r_{[Z]})$
inherits \textbf{lri}.
 We write  (\ref{autre1}) in the notation of
retractions and obtain
 the implications
\[ (\ref{autre1}) \Longleftrightarrow
[\alpha]^{[y]}=
[{\alpha}^y]=[\alpha]\Longleftrightarrow ^{lri} \;
[\alpha]=
{}^{[y]}{[\alpha]}\]
\[
 [\alpha]= {}^{[y]}{[\alpha]}\Longleftrightarrow^{Def.
\ref{[x]def}} \; {}^{\alpha}x ={}^{{}^y{\alpha}}x
\quad \text{for
all}\quad x \in X.\]
These equations together with part (1) imply
part (2) of the lemma.
\end{proof}

For instance, we have seen that $\Lcal_b,\Lcal_c$ in our example above are automorphisms. For our next example they are joined by $\Lcal_a=\Lcal_b\circ\Lcal_c$.

\begin{figure}
\begin{center}
\includegraphics[scale=.8]{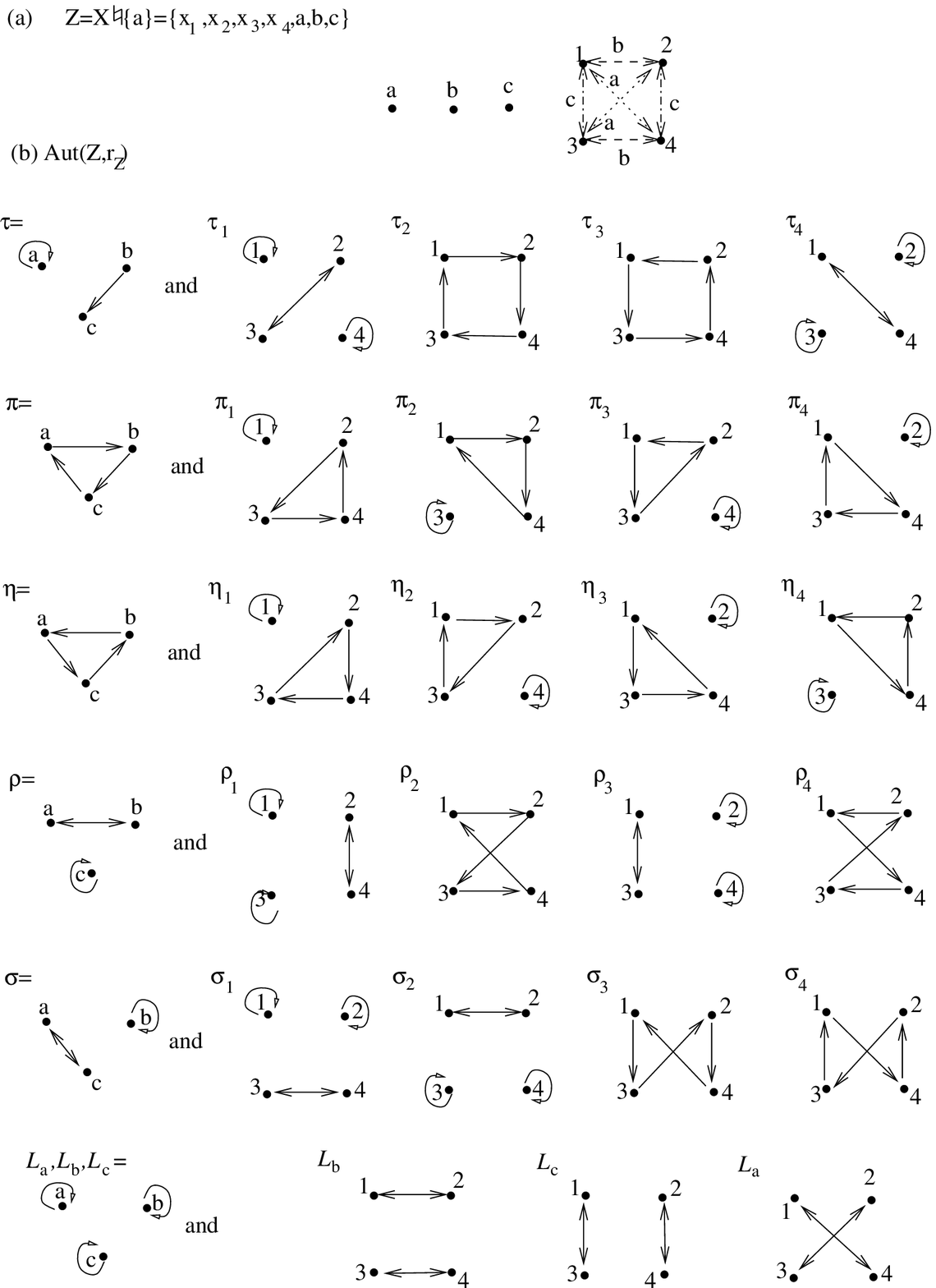}
\caption{(a) Graph of the solution of Example~\ref{autex2} and (b) its automorphism group.}
\label{figautex2}
\end{center}
\end{figure}

 \begin{example}
\label{autex2}  Let $(Z,r_Z)$ be the nondegenerate involutive square-free solution given by
\[
Z = \{x_1, x_2, x_3, x_4, a,b,c \},\  \Lcal_a=(x_1x_4)(x_2x_3),\
\Lcal_b=(x_1x_2)(x_3x_4), \ \Lcal_c=(x_1x_3)(x_2x_4).
\]
$(Z,r_Z)$ is an extension of the preceding solution $(X,r)$ by the trivial
solution on the one element set $\{a\}.$ Furthermore,
$\Lcal_a=\Lcal_b\circ \Lcal_c,$ so $\Gcal(Z,r_Z)=
\Gcal(X,r)\approx\mathbb{Z}_2 \times \mathbb{Z}_2.$ More
sophisticated arguments in Section 5 will show  that the group
$Aut(Z,r_Z)$ is isomorphic to $S_4$, the symmetric group on 4 elements.
 Furthermore, each automorphism $\tau \in Aut(X,r)$ can be
extended uniquely to an automorphism in $Aut(Z,r_Z),$
 by $\tau(a)=a.$ These are the automorphisms denoted $\tau_i$, $\Lcal_a,\Lcal_b,\Lcal_c$ shown in  Figure~\ref{figautex2}. The remaining automorphisms and
the graph of the solution are also shown.
\end{example}

We conclude the section with some straightforward generalities.
 \begin{remark}
Let $(X,r_X), (Y,r_Y)$ be solutions. By definition every
homomorphism of solutions $\varphi: (X,r_X) \rightarrow (Y,r_Y)$
agrees with the defining relations $\Re(X,r_X)$ and $\Re(Y,r_Y)$  of
the related algebraic objects, so it can be extended to
\begin{enumerate}
\item  a semigroup homomorphism $\varphi_{S}: S(X,r_X) \rightarrow
S(Y,r_Y)$ of their Yang-Baxter semigroups;
\item
a group homomorphism $\varphi_{G}: G(X,r_X) \rightarrow G(Y,r_Y)$
of their Yang-Baxter groups.
\item an algebra homomorphism $\varphi_{\Acal}: \Acal(k,X,r_X)
\rightarrow \Acal(k,Y,r_Y)$ of their Yang-Baxter algebras.
\end{enumerate}
If furthermore, $X$ is embedded in $G(X,r)$ (respectively in
$S(X,r)$), each $r$-automorphism $\tau$ of the solution $(X,r)$ can
be extendeded to a group automorphism $\tau_{G}: G(X,r) \rightarrow
G(X,r)$, (respectively, to a semigroup automorphism $\tau_{S}:
S(X,r) \rightarrow S(X,r)$. In this case, we have an embedding
$\Aut(X, r) \hookrightarrow \Aut(G(X,r))$, (respectively, an
embedding $\Aut(X, r) \hookrightarrow \Aut(S(X,r))$.
\end{remark}

\section{  Multipermutation solutions}
In this section we shall consider only nondegenerate symmetric sets
$(X,r)$.
\begin{definition}
\label{[x]def} \cite{ESS}, 3.2. Let $(X,r)$ be a nondegenerate
symmetric set. An equivalence relation $\sim$ is defined on $X$ as
\[ x \sim y \quad \text{ \emph{iff}} \quad
\Lcal x = \Lcal y.\] In this case we also have $\Rcal_x = \Rcal_y,$
see \cite{ESS}.

We denote by  $[x]$ the equivalence class of $x\in X$, $[X]=
X/_{\sim}$ is the set of equivalence classes.
\end{definition}
It is shown in \cite{ESS} that the solution $r$ induces a canonical
map
\[r_{[X]}: [X]\times[X] \longrightarrow [X]\times[X],\]
which makes $([X], r_{[X]})$ a nondegenerate symmetric set, called the
\emph{retraction of} $(X,r)$, and denoted $Ret(X,r).$

For our purposes we need concrete expressions of the left and the
right actions on $[X]$ canonically induced by  the left and right
actions on $X.$

\begin{definition}
\label{actionson[X]def} Let $(X,r)$ be a nondegenerate symmetric
set. Then the left and the right actions of $X$ onto itself induce
naturally left and right actions on the retraction $[X],$ via
\begin{equation}
\label{retractactionseq1}
{}^{[\alpha]}{[x]}:= [{}^{\alpha}{x}]\quad [\alpha]^{[x]}:= [\alpha^x], \;\text{for all}\; \alpha, x \in X.
\end{equation}
\end{definition}

Note that in the following there is no need to  assume \textbf{lri},
since $\Lcal_x=\Lcal_y$ \emph{iff} $\Rcal_x = \Rcal_y$.
\begin{lemma}
\label{retractactionslemma0} Suppose $(X,r)$ is a nondegenerate symmetric set.
 Then the left and the right actions
(\ref{retractactionseq1}) are well defined.
\end{lemma}
\begin{proof}
We need to show that ${}^{[\alpha]}{[x]}$ does not depend on the
representatives of $[\alpha]$ and $[x].$ It will be enough to show
that \[\Lcal_y = \Lcal_x \Longrightarrow [{}^{\alpha}{y}] =
[{}^{\alpha}{x}],\quad\text{for all}\; \alpha\in X.\] So fix
$\alpha, x \in X,$ and assume $y  \in [x].$ Let $z \in X.$ We have
to verify:
\begin{equation}
\label{reteq1}
{}^{({}^{\alpha}{y})}z= {}^{({}^{\alpha}{x})}z.
\end{equation}
 By hypothesis $(X,r)$ is nondegenerate,
so there exists a $t \in X,$ with
\begin{equation}
\label{reteq2}
z= {}^{\alpha^x}t= {}^{\alpha^y}t
\end{equation}
(The right hand side equality  comes from $\Lcal_y = \Lcal_x,$ which
implies   $\alpha^x=\alpha^y$). Then we have
\[
{}{}^{({}^{\alpha}{y})}z= {}^{{}^{\alpha}{y}}{({}^{\alpha^y}t)}
=^{\textbf{l1}}{}^{\alpha}{({}^yt)}=^{[x]=[y]} {}^{\alpha}{({}^xt)}=
{}^{{}^{\alpha}{x}}{({}^{\alpha^x}t)}=^{(\ref{reteq2})}\;{}^{({}^{\alpha}{x})}z.
\]
\end{proof}
\begin{lemma}
\label{retractactionslemma1} Suppose $(X,r)$ is a nondegenerate symmetric set.
Then
\begin{enumerate}
\item The left and the right actions (\ref{retractactionseq1})
define (as usual ) a canonical map \[r_{[X]}: [X]\times[X]
\longrightarrow [X]\times[X],
\]
which makes $([X], r_{[X]})$ a nondegenerate symmetric set.

\item \label{cyclicret} $(X,r)\; \text{cyclic} \Longrightarrow([X],
r_{[X]})\;\text{cyclic}$.

\item
\label{lriret} $(X,r)\; \text{is}\; {\bf lri} \Longrightarrow([X],
r_{[X]})\;\text{is}\; {\bf lri}.$
\item
\label{sqfreeret} $ (X,r)\; \text{square-free} \Longrightarrow ([X],
r_{[X]}) \; \text{square-free}.$
\end{enumerate}
\end{lemma}
\begin{proof}
It is easy to see that the left and the right actions
(\ref{retractactionseq1}) on [X] inherit conditions {\bf l1, r1,
lr3}. Therefore $r_{[X]}$ obeys YBE. The implication in
(\ref{lriret}) follows easily from (\ref{retractactionseq1}). We
leave conditions (\ref{cyclicret}), (\ref{sqfreeret}) (and so forth) to the reader.
\end{proof}

\begin{definition}
\label{mpldef}
 \cite{ESS}. The solution $([X], [r])$ is called
\emph{retraction of $(X,r)$} and is also denoted $Ret(X,r).$ For all
integers $m\geq 1$ $Ret^m(X,r)$ is defined recursively as
$Ret^m(X,r)= Ret(Ret^{m-1}(X,r)).$ A nondegenerate symmetric set
$(X,r)$ is \emph{a multipermutation solution of level} $m$ if $m$ is
the minimal number such that $Ret^m(X,r)$ is finite of order $1$,
this will be denoted by $mpl(X,r)= m.$ By definition $(X,r)$ is
\emph{a multipermutation solution of level} $0$ \emph{iff} $X$ is a
one element set.
\end{definition}

\begin{example}
\label{Lyubashenko}  An involutive
\emph{permutation solution} $(X,r)$ is defined \cite{Dri} (attributed to Lyubashenko) as
\[
r(x,y)=(\sigma(y),\sigma^{-1}(x)),
\]
where $\sigma$ is a fixed permutation in  $Sym(X).$
In this case $(X,r)$ is a nondegenerate symmetric set with
\textbf{lri}. Clearly, $\Lcal_x=\sigma,$ for all $x \in X,$ thus
$Ret(X,r)$ is a one element set and $mpl(X,r)=1.$ The converse also holds \cite{ESS},
so  $mpl(X,r)=1$ \emph{iff}  $(X,r)$ is a permutation
solution.
\end{example}

 In particular, the trivial solution $(X,r)$ with $X$ of order $\geq 2$ has $mpl(X,r)=1.$

\begin{remark}
\label{rethomlemma} Let $(X,r)$ be a nondegenerate symmetric set.
Then
 there is a
surjective homomorphism of solutions
\[
\mu: (X,r) \longrightarrow ([X], [r]); \quad \mu(x) = [x],
\]
Each finite symmetric set  $(X,r)$ with {\bf lri} has a well defined
oriented graph  $\Gamma(X, r_X)$, see \cite{TSh} and \ref{defgraph}.
In this case $\mu$ induces a homomorphism of graphs
\[
\mu_{\Gamma}: \Gamma(X, r_X) \longrightarrow \Gamma([X], r_{[X]})
\]
The graph $\Gamma([X], r_{[X]})$ is \emph{a retraction} of
$\Gamma(X, r_X).$
\end{remark}

The following  is straightforward from Definition \ref{mpldef}
\begin{lemma}
\label{mplretractlemma1} Suppose $(X,r)$ is a
 multipermutation solution. The following implications hold.
\[
 mpl(X,r) = m \;\Longrightarrow \;
mpl(Ret^k(X,r))= m-k,\quad \text{for all}\quad 1 \leq k \leq m-1.
\]
Conversely, if $k\geq 0$ is an integer, then
\[
mpl(Ret^k(X,r))= s  \;\Longrightarrow \; mpl(X,r) = s+k .
\]
\end{lemma}
\begin{lemma}
\label{retlemma1} Let $(X,r)$ be a nondegenerate square-free symmetric set. Then
 $[x]=[y], x,y \in X \Longrightarrow r(x,y)=(y,x).$
 \end{lemma}
 \begin{proof}
 Indeed,  $[x]=[y]$ implies ${}^yx = {}^xx = x,$ from which by  \textbf{lri},
 one has $x^y =  x$ and by an analogous argument with $x,y$ swapped one has $y^x=y$.
\end{proof}

\begin{corollary}
\label{mpl1lem} Suppose $(X,r)$ is a  nondegenerate square-free symmetric set (of order
$\geq 2$). Then the following conditions are equivalent.
\begin{enumerate}
\item
$mpl(X,r)=1$.
\item
$(X,r)$ is the trivial solution.
\item
${}^xy = y,$ for all  $x,y \in X.$
\item
$S(X,r)$ is the free abelian monoid generated by $X$.
\item
$G(X,r)$ is the free abelian group generated by $X$.
\item
$\Gcal(X,r) = \{ id_X \}$.
\end{enumerate}
\end{corollary}

 \emph{The $k$th retract orbit} of an element
 was introduced first
 in the case of square-free solutions, see \cite{T04}.

\begin{definition}  \label{kthretractdef}
   \cite{T04}
   Let $(X,r)$ be a nondegenerate symmetric set. We denote by $[x]_k$ the image of
$x$ in $Ret^{k}(X,r).$ The set
\[
\Ocal(x,k) := \{\xi\in X\mid [\xi]_k =[x]_k \}
\]
is called \emph{the $k$th retract orbit of} $x.$

Suppose  $X$ is an $r$-invariant subset of the solution $(Z,r)$, and
$x\in X.$ Then $[x]_{k,X}$ will denote the $k$-th retract of $x$ in
$X$.
\end{definition}

\begin{remark}
cf. \label{kthretractrem} \cite[Lemma 8.9]{T04} Let $(X,r)$ be a
square-free symmetric set.  For every positive integer $k\leq
mpl(X,r)$ the $k$th retract orbit $\Ocal(x,k)$ is $r$-invariant.
Furthermore if we denote by $r_{x,k}$ the corresponding restriction
of  $r$, then $(\Ocal (x,k), r_{x,k})$ is a multipermutation
solution  and
\begin{equation}
\label{mplretractorbiteq} mpl(\Ocal (x,k), r_{x,k})\leq k.
\end{equation}
\end{remark}
In \cite{T04}
(\ref{mplretractorbiteq}) is actually written as an equality, but  this should be
corrected as the following example shows.

\begin{example}
Let $X=\{x_1, \cdots, x_k, y\},$ and $r$ be defined via the actions
$\Lcal_y= (x_1 \cdots x_k), \Lcal_{x_i}=id_X.$ Then
$\Ocal(x_1,1)=\{x_1, \cdots, x_k\}, \Ocal(y,1)=\{y \}.$ Clearly,
$[X]=\{[x_1], [y]\}$ is the trivial solution, so $mplX =2.$ Note
that $mpl\Ocal(y,1)=0,$ $mpl(\Ocal(x_1,1))=1.$
\end{example}

\section{Strong twisted unions of solutions}

In this section we study special extensions of solutions called
\emph{strong twisted unions}. We recall first some basic facts and
definitions.

The notion of \emph{a union} of solutions and one-sided extensions
were introduced in \cite{ESS}, but only for nondegenerate involutive
solutions $(X,r_X),$ $(Y,r_Y).$ In \cite{TSh} are introduced and
studied more general extensions $(Z,r)$ of arbitrary solutions
$(X,r_X),$ $(Y,r_Y),$ and given necessary and sufficient conditions
(in terms of left and right actions) so that a regular extension
$(Z,r)$ satisfies YBE.

\begin{definition}
\label{extensiongeneraldef} \cite{TSh} Let $(X,r_X)$ and $(Y,r_Y)$
be disjoint quadratic sets (i.e. with bijective maps $r_X: X\times X
\longrightarrow X\times X, \; r_Y: Y\times Y\longrightarrow Y\times
Y$). Let $(Z,r)$ be a set with a bijection $r: Z\times
Z\longrightarrow Z\times Z.$ We say that $(Z,r)$ is \emph{a
(general)  extension of} $(X,r_X),(Y,r_Y),$ if $Z= X\bigcup Y$ as
sets, and  $r$ extends the maps $r_X$ and $r_Y,$ i.e. $r_{\mid X^2}=
r_X $, and $r_{\mid Y^2}=r_Y.$ Clearly in this case $X, Y$ are
$r$-invariant subsets of $Z$. $(Z,r)$ is \emph{a YB-extension of }
$(X,r_X)$, $(Y,r_Y)$ if $r$ obeys YBE.
\end{definition}

\begin{remark}
\label{extensionsrem} In the assumption of the above definition,
suppose $(Z,r)$ is a non-degenerate  extension of $(X,r_X),(Y,r_Y),$
(without any further restrictions on the solutions).  Then the
equalities $r(x,y) = ({}^xy,x^y),$ $r(y,x) = ({}^yx,y^x),$ and the
non-degeneracy of $r$, $r_X,$ $r_Y$ imply that
\[
{}^yx, x^y \in X, \quad {}^xy, y^x \in Y,\;\;
\text{for all}\;\; x \in X, y\in Y.
\]
Therefore, $r$ induces bijective maps
\begin{equation}
\label{rhosigma} \rho: Y\times X \longrightarrow X\times Y,\quad \sigma: X\times Y \longrightarrow Y\times X,
\end{equation}
and left and right ``actions"
\begin{equation}
\label{ractions1} {}^{(\;)}{\bullet}: Y\times X \longrightarrow
X,\;\;\; {\bullet}^{(\;)}: Y\times X \longrightarrow Y,\;
\text{projected from}\; \rho
\end{equation}
\begin{equation}
\label{ractions2}
 \la:
X\times Y \longrightarrow Y,\quad  \ra: X\times Y \longrightarrow X,
\ \text{projected\ from}\; \sigma.
\end{equation}
Clearly,  the 4-tuple of maps $(r_X, r_Y, \rho, \sigma)$ uniquely
determine the extension $r.$ The map $r$ is also uniquely determined
by $r_X$, $r_Y$, and the maps (\ref{ractions1}), (\ref{ractions2}).
\end{remark}

In the present paper  we restrict our attention to particular extensions  called
\emph{strong twisted unions}, also introduced in \cite{TSh}. However,  in the present
paper we prefer to avoid the most general form of this and focus on the case
where the extension is nondegenerate as in the remark above, and involutive.

\begin{definition} \cite{TSh}
\label{STUdef} In the notation of Remark \ref{extensionsrem} a
nondegenerate involutive extension $(Z,r)$ is a \emph{strong twisted union} of
the quadratic sets $(X,r_X)$ and $(Y, r_Y)$ if
\begin{enumerate}
\item
\label{a)}
The assignment $\alpha
\longrightarrow {}^{\alpha}\bullet $
extends to a left action of the
associated  group $G(Y, r_Y)$ (and the associated semigroup $S(Y,r_Y)$) on $X$,
and the assignment $x \longrightarrow {\bullet}^x $
extends to a
right action of the associated
group of $G(X,r_X)$  (and the associated semigroup $S(X,r_X)$)
  on $Y$;
\item
\label{b} The pair of ground actions satisfy
 \[\begin{array}{lclc}
 {\rm\bf stu :}\quad &{}^{{\alpha}^y}x = {}^{\alpha}x;
\quad  &{\alpha}^{{}^{\beta}x}={\alpha}^x,\quad  \text{for all}\quad x, y \in X, \alpha,\beta \in Y
\end{array}\]
\end{enumerate}
We shall use notation $(Z, r) = (X,r_X)\natural(Y, r_Y)$ (or
shortly, $(Z, r) =
      X\natural Y$) for a strong twisted union.
A strong twisted union  $(Z,r)$ of $(X,r_X)$ and $(Y, r_Y)$ is
 \emph{nontrivial} if at least one of the actions in
(\ref{a)}) is nontrivial.  In the case when
      both actions (\ref{a)}) are trivial we write $(Z, r) =
      X\natural_0 Y$). In this case one has
      $r(x,\alpha) = (\alpha,x)$ for all $x \in X,\alpha\in Y.$
 \end{definition}
The following example is extracted from \cite[Definition 3.3,
Proposition 3.9]{ESS},
 \begin{example}
\label{twisteduniondef}
Let  $(X,r_X), (Y,r_Y)$ be nondegenerate symmetric sets and $\sigma \in Aut(X,r), \rho \in Aut(Y,r).$
Define the nondegenerate involutive extension $(Z,r)$ via  $r_X,r_Y$ and the formulae
\begin{equation*}
r(\alpha, x)= (\sigma(x),\rho^{-1}(\alpha)),\quad
r(x,\alpha)=(\rho(\alpha),\sigma^{-1}(x)).
\end{equation*}
This is, moreover, a symmetric set (obeys the YBE) and is called a {\em twisted union} in \cite{ESS}. We denote it   $(Z,r) = X \natural_0 Y$ as a special case of a strong twisted union. Note that \[
\Lcal_{x\mid Z}=\Lcal_{x\mid X}.\rho, \quad \Lcal_{\alpha\mid
Z}=\Lcal_{\alpha\mid Y}.\sigma, \quad\; \text{for all}\; x\in X,
\alpha \in Y
\]
from which it is clear that $(Z,r)$ is nondegenerate as $(X,r_X)$ and $(Y,r_Y)$ are.
\end{example}
Clearly a  trivial extension $(Z,r)$ of $(X,r_X), (Y,r_Y)$ is a
particular case of twisted union. Another easy example is:

\begin{example} Let  $(X,r)$ be a nondegenerate symmetric set and
$Y= \{a\}$ a  one element  set with trivial solution. A strong twisted union (in fact any regular extension) of these that obeys the YBE is necessarily  a twisted union $Z=X\natural_0Y$ given by $\Lcal_a\in Aut(X,r)$ and  $id_Y$. This will be clear from  Proposition~\ref{stuybe} below. Our previous Example~\ref{autex2} is of this form.
\end{example}

The following lemma is straightforward.

\begin{lemma}
Suppose $(Z,r)=X \natural_0 Y,$ is a solution and $mpl(X,r_X)<
\infty, mpl(Y,r_Y) < \infty. $ Then $mpl(Z,r) = max\{mpl(X,r_X), \;
mpl(Y,r_Y)\}.$
\end{lemma}
In \cite[Definition 3.3]{ESS} the notion of \emph{a generalized
twisted union} $(Z,r)$ of the solutions $(X,r_X)$ and $(Y, r_Y),$ is
introduced in the class of symmetric sets.

\begin{definition}
\label{gtudef} A symmetric set $(Z,r)$ is \emph{a generalized
twisted union} of the disjoint symmetric sets  $(X,r_X)$, and $(Y,
r_Y)$ if it is an extension, and for every $x\in X, \alpha\in Y$ the
ground action ${}^{{\alpha}^x}{\bullet} : Y\times X \longrightarrow
X $ does not depend on $x$, and the ground action
${\bullet}^{{}^{\alpha}{x}} : Y\times X \longrightarrow Y$ does not
depend on  $\alpha.$
\end{definition}

\begin{remark}
\label{gtusturemark} Note that a strong twisted union $(Z,r)$ of
$(X,r_X)$ and $(Y, r_Y)$ does not necessarily obey YBE,
in contrast with twisted unions and  generalized twisted unions. It follows
straightforwardly from Definition \ref{gtudef} that a strong twisted
union $(Z,r)$ which is symmetric set  is a generalized twisted union
of $(X,r_X)$ and $(Y, r_Y)$.

Furthermore, it is shown in  \cite[ Proposition 8.3]{T04} that an
YB-extension $(Z,r)$ of two involutive square-free solutions
$(X,r_X)$, $(Y,r_Y)$ is a generalized twisted union   \emph{iff} it
is a strong twisted union. Proposition \ref{gtustuprop} now generalizes
this result for arbitrary symmetric sets with \textbf{lri}, without
necessarily assuming that the solutions are square-free.
\end{remark}

\begin{proposition}
\label{gtustuprop} Suppose the symmetric set  $(Z,r)$ has {\bf lri} and is an
extension of the solutions $(X,r_X)$, $(Y,r_Y)$. Then $(Z,r)$ is a generalized twisted union
\emph{iff} it is a strong twisted union.
\end{proposition}
\begin{proof}
By hypothesis $(Z,r)$ is an involutive solution with \textbf{lri},
thus Proposition \ref{lri&INVOLUTIVE} implies that $(Z,r)$ also
satisfies the cyclic conditions,  see Definition
\ref{cyclicconditionsall}. Assume $(Z,r)$ is a generalized twisted
union. Then by Definition \ref{gtudef} for every $x\in X, \alpha\in
Y$ the ground action ${}^{{\alpha}^x}{\bullet} : Y\times X
\longrightarrow X $ does not depend on $x$, and the ground action
${\bullet}^{{}^{\alpha}{x}} : Y\times X \longrightarrow Y$ does not
depend on  $\alpha.$ It follows then that for all $x,y \in X,
\alpha,\beta \in Y$ there are equalities:
\[
{}^{{\alpha}^x}{y}={}^{{\alpha}^y}{y}=^{\textbf{cl1}}\;
{}^{\alpha}{y} \quad \text{and} \quad {\beta}^{{}^{\alpha}{x}}=
{\beta}^{{}^{\beta}{x}}=^{\textbf{cr1}}\;{\beta}^x.
\]
We have shown that conditions \textbf{stu} are satisfied, it follows
then that $(Z,r) = (X,r_X)\natural(Y,r_Y).$ The converse implication
 is straightforward, see Remark
\ref{gtusturemark}.
\end{proof}

 Let $Ext^{\natural}(X,Y)$ denote the set of strong twisted unions which obey the YBE

 \begin{proposition}{\ }\label{stuybe}

\begin{enumerate}
\item
 A strong twisted union  of solutions
$(Z,r)=(X,r_X)\natural(Y,
r_Y)$ obeys YBE \emph{iff}

(a)  The assignment $\alpha \longrightarrow
{}^{\alpha}\bullet $
extends to a a group homomorphism

\[\varphi: G(Y,r_Y) \longrightarrow Aut(X,r);\quad
\text{and}\]

(b) The assignment $x \longrightarrow {\bullet}^x$
extends to a a
group homomorphism
\[ \psi: G(X,r_X) \longrightarrow Aut(Y,r).\]

In this case the pair of group homomorphisms
$(\varphi,\psi)$ is
uniquely determined by the pair of the ground actions,
or
equivalently by $r$.
\item
Furthermore, there is a one-to-one correspondence
between the sets
$Ext^{\natural}(X,Y)$ and $Hom(G(Y, r_Y),
Aut(X,r_X))\times Hom(G(X,
r_X), Aut(Y,r_Y))$  using (1).
\end{enumerate}
\end{proposition}

\begin{proof} We use \cite[Theorem~4.9]{TSh} which breaks down the
condition for any regular extension $(Z,r)$ (of which a strong twisted union is an example)
to obey the YBE into explicit conditions {\bf ml1,ml2,mr1,mr2}.
Particularly,
as  equalities in $ X \times X$
\[
{\bf ml2}:\quad {}^{\alpha}{r(x,y)}= r({}^{\alpha}{(x,y)})\quad
{\rm for\ all}\ x,y \in X, \alpha \in Y
\]
(the condition {\bf mr2} similar but for the right
action and the roles of $X,Y$ swapped.) We can
interpret this condition by computing both sides as
\[ {}^{\alpha}{r(x,y)}
={}^{\alpha}{({}^xy,x^y)}=({}^{\alpha}{{}^xy},{}^{(\alpha^{{}^xy}}{(x^y))}=
({}^{\alpha}{({}^xy)},{}^{\alpha}{(x^y))}=\Lcal_{\alpha}\times
\Lcal_{\alpha}\circ r(x,y)\] and
\[
r({}^{\alpha}{(x,y)})=r({}^{\alpha}x, {}^{{\alpha}^x}y)=
r({}^{\alpha}x, {}^{\alpha}y) =r\circ(\Lcal_{\alpha}\times
\Lcal_{\alpha}) (x,y).
\]
using our assumption {\bf stu}. Thus the condition has the meaning
\[
\Lcal_{\alpha}\times \Lcal_{\alpha}\circ r =
r\circ(\Lcal_{\alpha}\times \Lcal_{\alpha})
\]
that is $\Lcal_{\alpha\mid X} \in Aut(X,r),$ for every
$\alpha \in
Y$ as in part (a). Similarly for {\bf mr2} as in part (b). By definition every strong twisted
union obeys
{\bf ml1, mr1} so the above are the only conditions for $(Z,r)$ obeys YBE. Note that this proof works for general strong twisted unions as defined in \cite{TSh}  not only the nondegenerate involutive case in Definition~\ref{STUdef}, as long as the same definition is used for $Ext^\natural$.
\end{proof}


\begin{lemma}
\label{retractactionslemma2} Let $(Z,r)$ be a strong
twisted union of $(X,r_X)$ and $(Y, r_Y)$ which is a symmetric set with {\bf lri}.
Then the retraction $([Z], r_{[Z]})$ is a strong twisted union of
the retractions $([X],r_{[X]})$ and $([Y],r_{[Y]}).$
\end{lemma}
\begin{proof}
By hypothesis $Z$ is a  disjoint union of $X, Y$ thus [Z] is a
disjoint union of $[X], [Y].$ It follows from Lemma
\ref{retractactionslemma1} and conditions {\bf stu} on $Z$ that
\[
{}^{[\alpha]^{[y]}}{[x]}=^{(\ref{retractactionseq1})} [{}^{\alpha^y}{x}]
=^{{\bf stu}} [{}^{\alpha}x]=^{(\ref{retractactionseq1})}
{}^{[\alpha]}{[x]}.
\]
This gives the left hand side part of  {\bf stu}, which together
with {\bf lri} implies the right hand side of {\bf stu}. In view of Proposition~\ref{gtustuprop} this could be said equally well in terms of generalised twisted unions; the present version  can also be applied to general twisted unions obeying the YBE with {\bf lri} provided the retracts make sense. \end{proof}

In the following example the strong twisted union $Z=X\natural Y$
satisfies $mpl Z= max( mpl X, mpl Y).$
\begin{example}\label{examplestu1}
Let $Z=X \bigcup Y,$ where $(X, r_X)$, $(Y, r_Y)$, are the nondegenerate involutive square-free solutions defined by
follows.
$X= \{x,y, z\},$  $Y= \{\alpha,\beta,\gamma\},$ with
\[
\begin{array}{c}
\Lcal _{x\mid X}=(yz), \; \Lcal _{y\mid X}= \Lcal _{z\mid X}= id_X,\; \;
\Lcal _{\alpha\mid Y}= (\beta\gamma), \;\Lcal _{\beta\mid Y}=
\Lcal_{\gamma\mid Y}= id_Y.
\end{array}
\]
It is easy to see that $mpl X = mpl Y= 2.$ Define
\[
\begin{array}{c}
\Lcal _{x\mid Y}=(\beta\gamma),\; \Lcal _{y\mid Y}= \Lcal _{z\mid Y}=
id_Y, \; \; \Lcal _{\alpha\mid X}=(yz), \; \Lcal _{\beta\mid X}=
\Lcal_{\gamma\mid X}= id_X.
\end{array}
\]
and extend these to $Z$ (necessarily for any regular extension) by $ \Lcal_z=\Lcal_{z\mid X}\Lcal_{z\mid Y}$ where the restricted parts are considered to act trivially on the rest. So we have actions on $Z:$
\[
\begin{array}{c}
\Lcal _{x}=\Lcal _{\alpha}=(yz)(\beta\gamma),\; \;
\Lcal_{y}= \Lcal _{z}=\Lcal _{\beta}=
\Lcal_{\gamma}= id_Z,
\end{array}
\]
which define a nondegenerate involutive square-free solution $(Z,r).$ It is a strong
twisted union of $(X, r_X),$ and $(Y,r_Y).$ Clearly, $[x]=[\alpha],
[y]=[z]=[\beta]=[\gamma],$ thus $[Z]=\{[x],[y] \}$ is a two element
set, and $([Z],[r])$ is the trivial solution. This gives $mpl Z = 2
= mpl X = mpl Y.$

However, $(Z,r)$ can also be looked at as a strong twisted union of the
trivial solutions $X_0=\{x,\alpha\}, Y_0=\{y,z,\beta,\gamma\}$ with
\[
\begin{array}{lclc}
 \Lcal_{x\mid Y_0}=\Lcal_{\alpha\mid Y_0}=
(yz)(\beta\gamma);&\quad &\Lcal_{y\mid X_0}=\Lcal_{z\mid
Y_0}=\Lcal_{\beta\mid X_0}=\Lcal_{\gamma\mid X_0}= id_{X_0}.
\end{array}
\]
This way we have again $Z=X_0\natural Y_0,$ but
\[
mpl(Z,r)=mpl(X_0)+1=mpl(Y_0)+1.
\]
\end{example}

\begin{remark} \label{remarknonassoc}
Theorem \ref{mpl2theorem} in the next section shows that in the
case when $mpl(Z,r) = 2,$ $Z$ can always be presented as a strong
twisted union of a finite number of solutions $X_i, 1 \leq i \geq
s $ where $mpl(X_i), \leq 1,$  and there exists an $i$, such that
$mpl(X_i) \geq 1.$  The case $mpl(Z,r)\geq 3$ is more complicated.
We show in Proposition~\ref{mpl3prop} that every solution $(Z,r)$
with $mpl(Z,r)=3$ splits into $r$-invariant components $X_i,$
$1\leq i\geq s,$ where $X_i=\Vcal(\Gamma_i),$ $mpl X_i \leq 2,$
and for each pair $i,j, 1\leq i\leq s,$ the subset $X_{ij} =
X_i\bigcup X_j$ is $r$-invariant and has presentation as $X_{ij} =
X_i \natural X_j.$ However, in order to understand the nature of
solutions with higher multipermutation level and for their
classification we should understand under what conditions
$\natural$ associates or what replaces this. We believe that it is
always possible to present $(Z,r)$ as a strong twisted union of
components of strictly smaller multipermutation level.
\end{remark}

We now give two more examples of  strong twisted
unions of solutions. In the first we have a strong twisted union
$Z=X\natural Y,$ with $mpl X=2, mpl Y = 1$ and $mpl Z
= 3.$ $Z$ is also split as a strong twisted union of three
$r$-invariant components, each of multipermutation level 1, and the
``associative law" holds. The second example, see Example
\ref{exCUBE}, gives a solution $(Z,r)$ as an extension (but
not necessarily a strong twisted union) of two
$r$-invariant subsets. It illustrates Remark \ref{remarknonassoc}
with respect to $\natural$ not  being nonassociativite in general.

\begin{example}\label{examplestu2}
Let $(X,r_X)$  be the nondegenerate involutive
square-free solution
defined by
$X=\{a, b, c, x_1, x_2, y_1, y_2, z_1, z_2\}$ and
left action:
\begin{equation}
\label{exeq0} \Lcal_{a\mid X}=\Lcal_{b\mid
X}=\Lcal_{c\mid X}=(x_1x_2)(y_1y_2)(z_1z_2),
\end{equation}
and $ \Lcal_{x_i\mid X} =
\Lcal_{y_i\mid X} =
\Lcal_{z_i\mid X}=id_X$ for $i=1,2$.

Let $(Y, r_Y)$ be the trivial solution on the set
$Y=\{\alpha,\beta\}$.
We extend the left actions on $Z=X \bigcup
Y$ as
\[
\Lcal_{\alpha}=\Lcal_{\beta}=(abc)(x_1 y_1 z_1 x_2y_2
z_2),
\]
\[
\Lcal_{a\mid Y}=\Lcal_{b\mid Y}=\Lcal_{c\mid Y}=
\Lcal_{x_i\mid Y}
= \Lcal_{y_i\mid Y} = \Lcal_{z_i\mid Y}=id_Y, i=1,2.
\]
This defines a left action of  $Z$ which we verify
defines $(Z,r)$ as a nondegenerate
involutive square-free solution. Clearly, $mpl
Y= 1,$ and it is easy to see that $mpl(X, r_X)= 2$. We
leave the reader to verify that $mpl(Z, r)= 3$.

\begin{figure}\begin{center}
\includegraphics[scale=1]{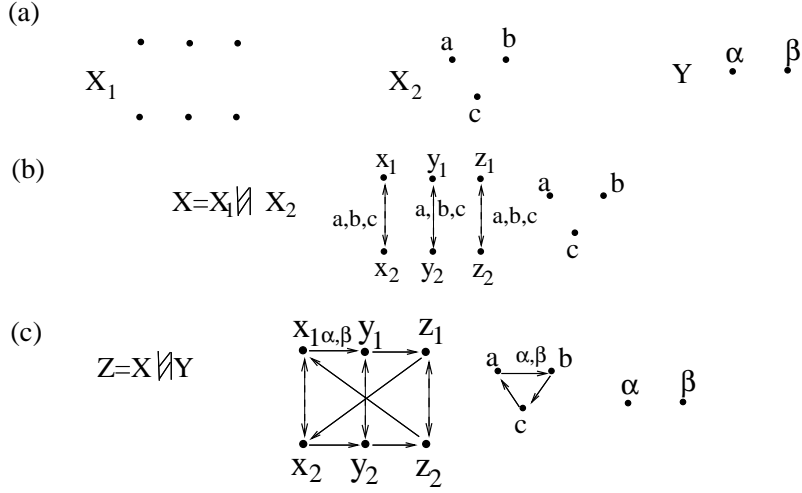}
\caption{Graph of Example~\ref{examplestu2} of a strong twisted union}
\label{figex3} \end{center}
\end{figure}

The graph $\Gamma(X, r_X)$ has three nontrivial
components and three one vertex components, as shown in Figure~\ref{figex3}
part (b). Consider the presentation of $X = X_1 \bigcup X_2$,
where $X_1=\{ x_1, x_2, y_1, y_2, z_1, z_2\}$,
$X_2=\{a, b, c\}$. Both are $r$-invariant sets and the
restrictions $r_1= r_{\mid X_1\times X_1}$, $r_2= r_{\mid
X_2\times X_2}$, are the trivial solutions as shown along with $Y$ in
part (a) of the figure. Moreover, the actions given by
(\ref{exeq0}) make $X$ a strong twisted union $X = X_1\natural X_2.$
This way we have
\[
Z = (X_1\natural X_2)\natural Y = X_1\natural
(X_2\natural Y) =X_2\natural (X_1\natural Y)\]
\[ mpl X_1 = mpl X_2 = mpl Y = 1, mpl(Z)= 3.\]
The graph $\Gamma(Z,r)$ is shown in part (c) of the
figure.
\end{example}

\begin{example}
\label{exCUBE} Let  $X_3 =\{a \} $ be the one element
solution ($r_3
= id_{X_3 \times X_3}$), and let $(X_1, r_1), (X_2,
r_2)$ be the
trivial solutions on the sets
\[
X_1 = \{x_1, x_2, x_3, x_4, x_5, x_6, x_7, x_8 \},
\quad X_2 =
\{a,b,c \}
\]
as depicted in Figure~\ref{figexcube} (a). Consider the strong twisted union $(X, r_0)= X_1
\natural X_2$  with $r_0$
defined by the left actions   \[
\Lcal_b=(x_1x_2)(x_3x_4)(x_5x_6)(x_7x_8),\quad
\Lcal_c=(x_1x_5)(x_2x_6)(x_3x_7)(x_4x_8)
\]
and the condition {\bf lri}. As usual, in this way we find a nondegenerate involutive square-free solution, as depicted in part (b). It is easy to see that
$Ret(X,r_0)$ is the trivial solution on the
set $[X]=\{[x_1],[b],[c]\}$, therefore $mpl(X,r_0)=2.$

Consider now the set $Z = X \bigcup X_3.$ Let
\[
\Lcal_a=(bc)(x_1x_4)(x_2x_3)(x_5x_8)(x_6x_7).\
\]
We leave the reader to verify that $\Lcal_a \in Aut
(X,r_0).$ Next
we define the extension $(Z,r) = X \natural X_3,$ with
$r$ defined
via $\Lcal_a$, the condition {\bf lri} that we are extending $r_0, r_3.$
Clearly, $(Z,r)$ is nondegenerate
involutive square-free solution and
$mpl(Z,r)=3$. This way we have a presentation which we denote $Z=
(X_1\natural
X_2)\natural X_3= (X_2\natural X_1)\natural X_3$ as in part (c).

\begin{figure}\begin{center}
\includegraphics[scale=.9]{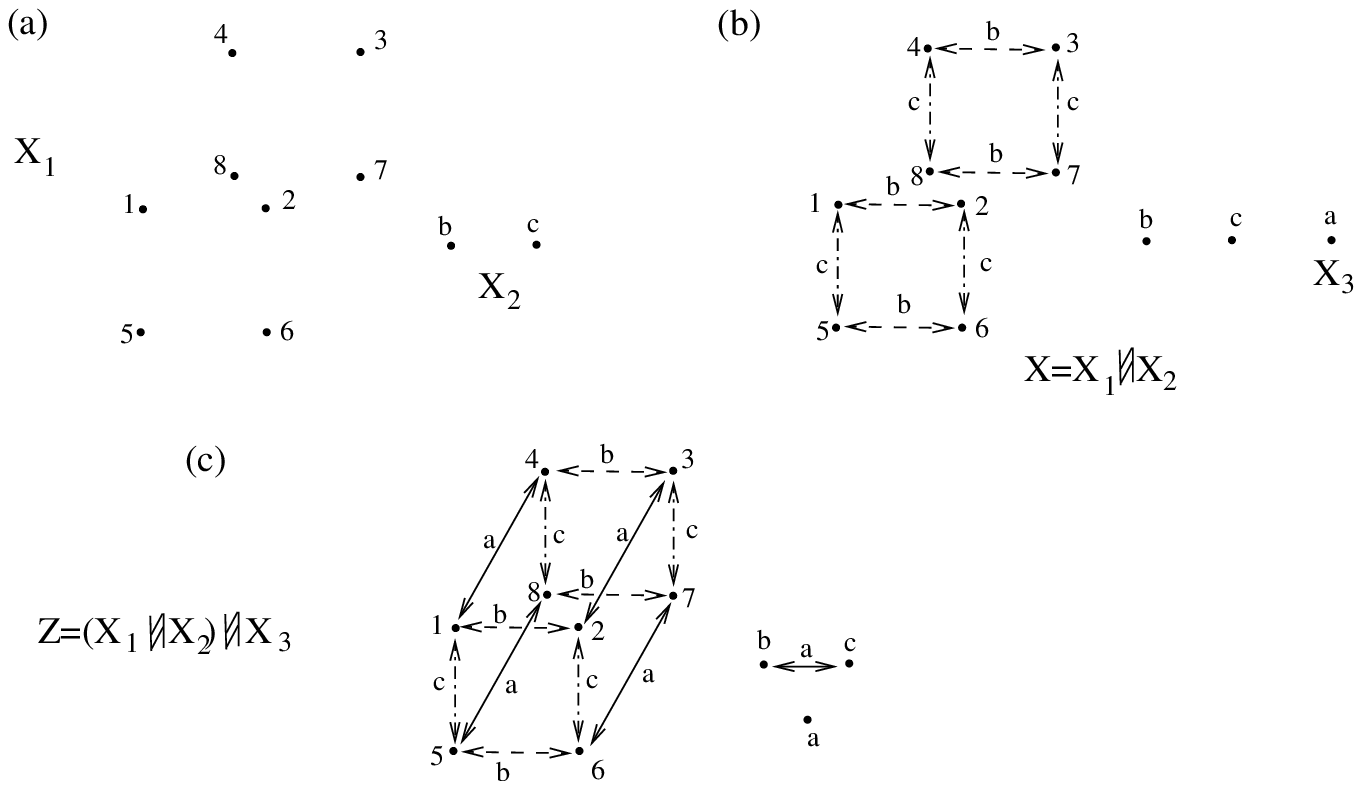}
\caption{Graph of Example~\ref{exCUBE}}
\label{figexcube} \end{center}
\end{figure}

Consider now  the  set $Y =X_1\bigcup X_3,$ which is
an
$r$-invariant subset of $(Z, r).$ Let $r_Y= r_{\mid Y
\times Y},$
then $(Y, r_Y) = X_1\natural X_3,$ ($r_Y$ can be
defined also via
the action $\Lcal_{a\mid X_1}=
(x_1x_4)(x_2x_3)(x_5x_8)(x_6x_7)$). One
has $mpl(Y,r_y) = 2.$ Clearly $(Z, r)$ as an extension of $(Y,r_Y)$ and
$(X_2, r_2).$
However, this extension is not a strong twisted union,
since neither
$\Lcal_b$ nor $\Lcal_c$ is in $Aut(Y, r_Y).$  Or
equivalently, a
direct verification shows that the pair $Y, X_3$ does
not satisfy
condition \textbf{stu}, since
\[
{}^{b^a}{x_1}= {}^{c}{x_1}= x_5 \neq {}^{b}{x_1}=x_2.
\]
Thus, while one also has $Z=X_1 \natural(X_3 \natural
X_2)$ by similar computations, the parentheses do not exactly associate
with $\natural$; we have shown that $(Z, r)$ is
a YB-extension of the nondegenerate square-free involutive
solutions $Y=X_1 \natural X_3$ and $X_2$ as desired but this extension is
not a strong twisted union.
\end{example}

\begin{theorem}
\label{mplstutheorem} Let  $(X,r_X)$ and $(Y, r_Y),$
be square-free multipermutation solutions, with $mpl X=m$, $mpl Y=k$.
 Suppose $(Z, r) = X\natural Y$ is a nondegenerate symmetric set. Then
\begin{equation}
\label{mplstueq} mpl(Z, r) \leq max\{m, k \}+1.
\end{equation}
\end{theorem}
\begin{proof}
Without loss of generality we can assume  $k \leq m.$
We shall prove (\ref{mplstueq}) by induction on $m.$ For the base for the
induction, $m=1$, we have to show that $mpl(Z,r) \leq 2,$ which  by
Lemma~\ref{mplstueq} is equivalent to $mpl(Ret(Z,r)) = 1.$
By hypothesis \textbf{stu} holds on $Z,$ so
\[
{}^{{\alpha}^y}x= {}^{\alpha}x,\quad
{}^{x^{\alpha}}{\beta}={}^{x}{\beta}, \; \text{for every} \; x,y \in X, \;  \alpha,\beta \in Y,
\]
or equivalently,
\begin{equation}
\label{stueq2} \Lcal_{{\alpha}^y\mid
X}=\Lcal_{\alpha\mid X}=\Lcal_{{}^y{\alpha}\mid X}; \;
\Lcal_{x^{\alpha}\mid Y}=\Lcal_{x\mid Y}=\Lcal_{{}^{\alpha}x\mid Y},\;
\text{for all} \;
x,y \in X, \alpha,\beta \in Y.
\end{equation}
By Corollary  \ref{mpl1lem} $mpl X =1 $ implies  $\Lcal_{x\mid X}=
id_X$.  As a set $Z$ is a disjoint union $Z =X\bigcup
Y,$ hence the following equalities hold.
\begin{equation}
\label{stueq3} \Lcal_{x\mid Z}=\Lcal_{x\mid
X}\circ\Lcal_{x \mid Y}= id_X\circ\Lcal_{x\mid Y}= \Lcal_{x\mid Y}, \quad
\forall x \in X.
\end{equation}
Analogously,
\begin{equation}
\label{stueq4} \Lcal_{\alpha\mid Z}=\Lcal_{\alpha\mid X},
\quad\forall \alpha \in Y.
\end{equation}
Now from (\ref{stueq2}), (\ref{stueq3}), and (\ref{stueq4}) we
deduce the equalities
\begin{equation}
\label{stueq5} \Lcal_{{\alpha}^y\mid Z}=\Lcal_{\alpha\mid
Z}=\Lcal_{{}^y{\alpha}\mid Z}; \quad\Lcal_{x^{\alpha}\mid
Z}=\Lcal_{x\mid Z}=\Lcal_{{}^{\alpha}x\mid Z},\;
\text{for all } \; x,y \in X, \alpha,\beta \in Y.
\end{equation}
This, together with  Corollary~\ref{mpl1lem} imply
that either i) $Ret(Z,r)$ is a trivial solution of order $\geq 2,$
hence $mpl(Z,r)=2$; or ii) $Ret(Z,r)$ is one element
solution so  $mpl(Z,r) = 1.$ We have verified (\ref{mplstueq}) for
$m =1.$ Assume now (\ref{mplstueq}) holds for all $m \leq m_0.$
Suppose $mpl X =m_0+1,$ thus by Lemma \ref{mplstueq},
$mpl([X], r_{[X]}) = m_0.$ Clearly, $mpl([Y], r_{[Y]})
= mpl(Y, r_Y)-1 \leq m_0.$ By Lemma \ref{retractactionslemma2} the
retract $([Z], r_{[Z]})$ is a strong twisted union $([Z], r_{[Z]})=
[X]\natural[Y]$, hence by the inductive assumption one
has $mpl([Z], r_{[Z]}) \leq m_0 +1.$ This and the equality
$mpl(Z,r)=mpl([Z], r_{[Z]})+1$ implies $mpl(Z,r)\leq m_0 +2.$ It
follows then that (\ref{mplstueq}) holds for all $m, m \geq 1,$
which proves the theorem.
\end{proof}

The following lemma is straightforward.
\begin{lemma}

\label{mplstulemma}
Let $(X_1,r_1), (X_2,r_2)$ be disjoint solutions with
$mpl(X_1,r_1)= m_1$, $mpl(X_2,r_2)= m_2$, let $m = max \{m_1, m_2\}$.
Suppose $(Z, r) =X_1\natural X_2$ is a symmetric set. Then $mpl Z = m +1$
\emph{iff} for some $i, 1 \leq i \leq
2$ with $m_i = m$ there are $x \in X_i,$ and  $\alpha
\in Z \backslash X_i$ such that the orbit $O^{m-1}(x)$ does not
contain ${}^{\alpha}x.$
\end{lemma}

Lemma \ref{mplstulemma} has a clear interpretation in terms of the
graphs $\Gamma(X_i, r_i), i =1,2, \Gamma(Z, r),$ see Definitions~\ref{defgraph},
 \ref{defGammax}.
In the case discussed by the lemma, $\Lcal_{\alpha\mid X_i}$ acts as
an automorphism of $(X_i, r_i)$ which maps the connected component
$\Gamma_x$ of $\Gamma(X_i, r_i)$ onto a different isomorphic component
$\Gamma_{{}^{\alpha}x}$ of $\Gamma(X_i, r_i).$
Clearly this is possible only in the case when $\Gamma(X_i, r_i)$
 has at least two connected components which are isomorphic as
 graphs. The last is a necessary but not a sufficient condition.

\section{Graphs of symmetric sets with {\bf lri} }

Each finite involutive solution $(X,r)$ with {\bf lri} can be
represented geometrically by its \emph{graph of the left action}
$\Gamma(X,r).$ It is an oriented labeled multi-graph (although we
refer to it as a {\em{graph}}). It was introduced  in \cite{T00} for
square-free solutions, see also \cite{TSh}. Here we recall the
definition.

\begin{definition}\cite{TSh} \label{defgraph}
Let $(X,r)$ be a finite symmetric set  with
 \textbf{lri}, we define the (complete)
graph $\Gamma=\Gamma (X,r)$ as follows. It is an oriented graph,
which reflects the left action of
 $G(X,r)$ on $X$. The set of
vertices of $\Gamma$ is exactly $X.$ There is a labeled arrow  $x
{\buildrel a \over \longrightarrow}y,$ if $x,y,a \in X$
and ${}^ax = y.$ An edge $x {\buildrel a \over \longrightarrow}y,$
with  $x\neq y$ is called \emph{a nontrivial edge}. We will often
consider the simplified graph in which to avoid clutter we typically
omit self-loops unless needed for clarity or contrast. Also for the
same reason, we use the line type to indicate when the same type of
element acts, rather than labeling every arrow. Clearly,
$x{\buildrel a\over \longleftrightarrow} y$ indicates that ${}^ax =
y$ and ${}^ay = x.$ (One can make such graphs for arbitrary solutions
but then it should be indicated which action is considered).
\end{definition}

Note that two solutions are isomorphic if and only if their
(complete) oriented graphs are isomorphic. Various properties of a
solution $(X,r)$ are reflected in the properties of its graph
$\Gamma(X,r)$, see for example the remark below,
Proposition~\ref{graphlemma1}, Theorem \ref{graphprop1}.

\begin{remark} Let $(Z,r)$ be a symmetric set with \textbf{lri},
 $\Gamma=\Gamma(Z,r)$.
\begin{enumerate}
\item $(Z,r)$ is a square-free solution \emph{iff}
$\Gamma$ does not contain a nontrivial edge $x{\buildrel x\over
\longrightarrow} y, x \neq y.$
\item In this case, $(Z,r)$ is a trivial solution, or equivalently, $mpl(Z,r)=1$
\emph{iff} $\Gamma$ does not contain nontrivial edges $x {\buildrel
a \over \longrightarrow}y, x \neq y$ .
\end{enumerate}
\end{remark}

By Proposition \ref{lri&INVOLUTIVE} our assumptions that $(X,r)$ is
a symmetric set with {\bf lri} imply the cyclic conditions (without
necessarily assuming $(X,r)$ square-free). We will need the
assumption $(X,r)$ square-free any time when we claim ${}^xx=x,$ for
every $x\in X.$ Examples of graphs of square-free solutions were already given in Section 2, see Examples
\ref{autex1}, and \ref{autex2}, as well as in Section~4.

We will find now various properties of $\Gamma= \Gamma(X,r).$
\begin{notation}
\label{defGammax}
Suppose $\Gamma_0$ is a subgraph of $\Gamma.$ Denote by
$\Vcal(\Gamma _0)$ the set of all vertices of $\Gamma_0.$
$\Ecal(\Gamma_0)$ denotes the set of all labels of (nontrivial)
edges that occur in $\Gamma_0,$ i.e.
\[
\Ecal(\Gamma_0)= \{ a \in X\mid \exists \; \text{an edge}\; x
{\buildrel a \over \longrightarrow}y \subset \Gamma _0, x \neq y\}.
\]
Clearly, each $x\in X$ determines uniquely a  connected component
 of  $\Gamma$ which contains $x$ as a vertex, we shall denote it by
 $\Gamma_x$.
\end{notation}

Let $a,x \in X.$ Suppose ${}^ax\neq x.$ Then the orbit of $x$ under
the left action of the cyclic group ${}_{gr} \langle a \rangle$ (or
equivalently under the action of $\Lcal_a$) on $X$ is a cycle $(x_1
x_2 \cdots x_m)$ of length $m \geq 2$ in the symmetric group
$Sym(X)$, where for symmetry we set $x_1 = x.$ One has ${}^ax_i =
x_{i+1}, 1 \leq i\leq m-1, {}^ax_m =x_1,$. This cycle participates
in the presentation of the permutation $\Lcal_a \in Sym(X)$ as a
product of disjoint cycles.
 Clearly, $x_2 \cdots x_m \in \Gamma _x.$

\begin{notation}
\label{Lax} For $x, a$ as above we use notation $\Lcal_a^{x}=
(x_1\cdots x_m) = \Lcal_a^{x_i}$, $1\leq i \leq m. $
\end{notation}

\begin{convention}
Till the end of the section  we shall consider only  finite
square-free nondegenerate symmetric sets  $(X, r)$. We recall from
Theorem \ref{basictheorem} that such solutions satisfy both
\textbf{lri} and the cyclic conditions. In particular, one has
\begin{equation}\label{remarkoncycles} \Lcal_{(a^x)}^x=
\Lcal_{a}^x = \Lcal_{({}^xa)}^{x}.\end{equation} Moreover, it is
known that for each such square-free solution $(X,r)$ with $1 <
|X|< \infty$, its YB  group $G(X,r)$ acts nontransitively on $X$,
see \cite{T04, Rump}, and therefore the graph $\Gamma(X,r)$ has at
least two connected components.
\end{convention}

\begin{proposition}
\label{graphlemma1} Let $(X,r)$ be a finite nondegenerate
involutive square-free solution, $G=G(X,r),$ be its YB group, and
$\Gamma=\Gamma(X,r)$ be its graph. Suppose  $\Gamma_1, \Gamma_2,
\cdots, \Gamma_s$ is the set of all connected components, with
sets of vertices, respectively, $X_i=\Vcal(\Gamma_i), 1 \leq i
\leq s.$
\begin{enumerate}
\item
Each set $X_i,  1 \leq i \leq s,$ is precisely an orbit of the
left action of the group $G$ on $X$.
\item
 $X_i$ is $r$-invariant, so $(X_i, r_i)$, with $r_i = r_{\mid X_i\times X_i}$,
 is a nondegenerate involutive square-free solution. Its graph $\Gamma(X_i, r_i)$  is the
subgraph of $\Gamma_i$ obtained by erasing the edges labeled
with elements $a \in \Ecal(\Gamma_1) \backslash X_i.$
\item Furthermore,
if $mpl (X,r) = m,$ then $mpl(X_i, r_i) \leq m-1,$ for each $i, 1
\leq i \leq s.$ More precisely,  in this case,  $x \in X_i$ implies
$X_i \subseteq \Ocal(x, m-1).$
\end{enumerate}
\end{proposition}
\begin{proof}
Sketch of the proof. Condition (1)  follows straightforwardly from
the definition of the graph $\Gamma(X,r)$. The definition of
$\Gamma(X,r)$ also implies that the set of vertices $X_i$ of each
connected component is invariant under the left action. By
\textbf{lri} it is also invariant with respect to the right
action, and therefore it is $r$-invariant. We shall prove (3). By
\ref{mplretractlemma1} one has $mpl(Ret^{m-1}(X,r))= 1.$
$Ret^{m-1}(X,r)$ inherits the property of being square-free, so by
Lemma \ref{mpl1lem}, $Ret^{m-1}(X,r)$ is a trivial solution. In
terms of the left action this is equivalent to
\[
{}^yx \in \Ocal(x, m-1)\quad \text{for all}\quad x,y \in X.
\]
It follows then that the orbit of $x,$ $X_i= \Ocal_G(x)$ is
contained in the $(m-1)$st retract orbit $\Ocal(x, m-1),$ and Remark
\ref{kthretractrem}, implies $mpl(X_i)\leq m-1.$
\end{proof}

\begin{definition}
 \label{locallycommutingactionsdef}
 Let $a\in
X,$ and let  $\Gamma_1$ be a connected component of  $\Gamma$. We
shall use notation $\Lcal_{a\mid \Gamma_1}$ for the restriction of
the action $\Lcal_a$ on the set of vertices $\Vcal(\Gamma_1).$ We
say that \emph{the (left) actions of $a$ and $b$ commute on
$\Gamma_1$} if the following equalities hold:
\begin{equation}
\label{commutingeq}
 \Lcal_{({}^ba)\mid \Gamma_1} =\Lcal_{a\mid \Gamma_1}, \quad
 \Lcal_{({}^ab)\mid \Gamma_1} =\Lcal_{b\mid\Gamma_1}.
\end{equation}
In this case we also say that each two  edges of  $\Gamma_1$ labeled
by $a$ and $b$ commute.
\end{definition}

The graph $\Gamma(X,r)$ and the combinatorial properties of the
solution $(X,r)$ can be used as a ``mini-computer" to deduce missing
relation of $(X,r)$. For example, see the proofs of Lemmas
\ref{graphcellslemma}, \ref{graphlemma2}, \ref{graphlemma3}
and Example \ref{examplempl3}.

The following straightforward lemma describes the ``cells" from which
we build $\Gamma(X,r)$ in the general case of $(X,r)$ with
\textbf{lri}.
\begin{lemma}
\label{graphcellslemma} Let $(X,r)$ be a symmetric set with
\textbf{lri}.
\begin{enumerate}
\item
For each $x_1, a, b \in X, a \neq b,$  there exists a uniquely
determined subgraph $\Gamma _0$ of $\Gamma$:
\begin{equation}\label{ecprop1lemadefaut1}
\xymatrix{
x_1 \ar[r]^a \ar[d]_ {b} & x_2 \ar[d]^{{}^ab} \\
x_3 \ar[r]_{{}^ba} & x_4 },
\end{equation}
where the vertices $ x_1, x_2, x_3, x_4 \in X$ are not necessarily
pairwise distinct.
\item
Given $a, b,$ each vertex $x_i, $ of $\Gamma _0$ in
(\ref{ecprop1lemadefaut1}) determines the remaining three verteces
$x_j, 1\leq j\leq 4, j \neq i.$
\item
\label{graphcells3} In the notation of (1) let $\Gamma_{x_1}$ be the
connected  component which contains $x_1$ as a vertex. Suppose that
the left actions of $a$ and $b$ commute on $\Gamma_{x_1}$, see
Definition \ref{locallycommutingactionsdef}. Then the subgraph
$\Gamma _0$ has the shape
\begin{equation}\label{ecprop1lemadefaut2}
\xymatrix{
x_1 \ar[r]^a \ar[d]_ {b} & x_2 \ar[d]^{b} \\
x_3 \ar[r]_{a} & x_4 },
\end{equation}
and the following implications hold:

\begin{equation}
\label{ecprop1lemadefaut3}
\begin{array}{lclc}
x_1\neq x_2 \quad & \Longleftrightarrow \quad& x_3\neq x_4\\
x_1\neq x_3 \quad & \Longleftrightarrow \quad& x_2\neq x_4.
\end{array}
\end{equation}
\end{enumerate}
\end{lemma}
\begin{proof}
(1) is clear. Note that given $\xi, a \in X $ there exists unique
$x \in X,$ such that ${}^ax=\xi.$ Indeed, by \textbf{lri}, $x
={\xi}^a$. One uses this and \textbf{lri} again to deduce  (2).
Now we prove the first implication of (3). It follows from the
diagram (\ref{ecprop1lemadefaut2}) that one has $x_3={}^b{x_1},
x_4={}^{{}^ab}{x_2}={}^{b}{x_2}$, so by \textbf{lri}
${x_3}^b=x_1$, ${x_4}^b=x_2.$ Clearly this implies
\[
x_1 = x_2 \quad \Longleftrightarrow\quad  x_3 = x_4.
\]
Analogous argument proves the second implication.
\end{proof}

Note that the diagram (\ref{ecprop1lemadefaut1}) is just a formal
graphic expression of condition \textbf{l1}, so it always presents
a subgraph of $\Gamma(X,r),$ although in the case when some
vertices coincide it is deformed to a segment or to a single vertex
with loops. It is convenient (and mathematically correct) for
computational purposes to use sometimes diagrams of the shape of
squares (as in (\ref{ecprop1lemadefaut1})) even when some vertices
coincide. See for example Corollary \ref{graphcor1}.

\begin{example}
\label{examplempl3} Let $(Z,r)$ be a nondegenerate involutive square-free solution with at
least $7$ elements. Suppose  where $a,b,c, x_1, x_2, x_3 \in Z,$ and
\begin{equation}
\label{exeq1}
\begin{array}{lclc}
r(a,b)=(c,a),\quad & r(a,c)=(b,a), \quad& r(a,x_1)=(x_2,a)\quad \\
r(a,x_2)=(x_1,a),\quad & r(a, x_3)= (x_3,a),\quad & r(b,c)=(c,b),\quad
r(b,x_1)=(x_3,b)
\end{array}\end{equation}
Then $r(c,x_2)=(x_3,c)$ and there exists some $x_4\in X,$ distinct from
$a,b,c, x_1, x_2, x_3,$ such that
\begin{equation}
\label{exeq2}
\begin{array}{lclc}
r(c,x_1)=(x_4,c), \quad & r(b,x_2)=(x_4,b),\quad& r(a,x_4)=(x_4,a).
\end{array}\end{equation}
This is indicated in terms of the graph in Figure~\ref{figex4}.
\begin{figure}
\begin{center}
\includegraphics{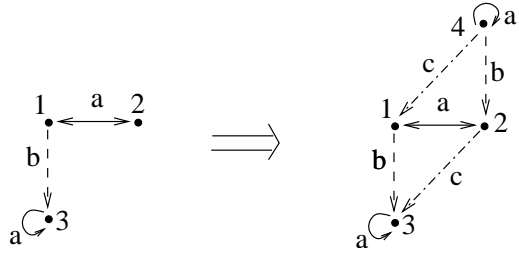}
\caption{Diagram for Example~\ref{examplempl3}}
\label{figex4}\end{center}
\end{figure}
Furthermore, if $(Z,r)$ has exactly $7$ elements then
(\ref{exeq1}) determines uniquely the solution $(Z,r).$ It is not
difficult to see that in addition to  (\ref{exeq1}), (\ref{exeq2}),
$(Z,r)$ satisfies $r(x_i, x_j)=(x_j,x_i),$ for all $1\leq i,j \leq
4.$ For the left actions we have:
\begin{equation*}
\begin{array}{lclc}
 \Lcal_{x_i}= id_Z, 1\leq i\leq 4,\quad &  \Lcal_{b}=(x_1x_3)(x_2x_4),
 \quad \Lcal_{c}=(x_1x_4)(x_2x_3),\quad &\Lcal_{a}=(bc)(x_1x_2).
\end{array}\end{equation*}
So the retracts of $(Z,r)$ are:
\begin{equation*}
\begin{array}{lclc}
Ret(Z,r) =\{[x_1], [a], [b],
[c]\},\quad&\Lcal_{[a]}=([b][c]),\quad&\Lcal_{[b]}=\Lcal_{[c]}=\Lcal_{[x_1]}=id_{[Z]}
 \\Ret^2(Z,r) =\{[a]^{(2)}, [b]^{(2)}\},\quad& \text{the trivial solution}&
 Ret^3(Z,r) =\{[a]^{(3)}\}.
\end{array}\end{equation*}
It follows then that $mpl(Z,r)=3.$

\end{example}
\begin{example}
\label{examplempl4} The solution from the previous example can be
constructed as a strong twisted union of disjoint solutions of lower
multipermutation levels. Consider the trivial solutions $(X_1,
r_1), (X_2, r_2),$ where
\[
X_1 =\{ x_1, x_2, x_3, x_4 \}, \quad X_2 =\{ b, c \}.
\]
Let $(X,r_X)=X_1\natural X_2$ be the strong twisted union defined
via the actions
\begin{equation*}
\begin{array}{lclc}
\Lcal_{b}=(x_1x_3)(x_2x_4),
\quad&\Lcal_{c}=(x_1x_4)(x_2x_3),\quad&\Lcal_{x_i}=id_X.
\end{array}\end{equation*}
We leave the reader to verify that $(X,r_X)$ is a solution.
 The graph $\Gamma(X,r_X)$ has three connected
components: two one vertex components, $b$ and $c$ and $\Gamma_1,$
with $\Vcal(\Gamma_1)=X_1$, $ \Ecal(\Gamma_1)=X_2.$
Moreover, $\Gamma_1$ is a `graph of first type', see Definition
\ref{graph1} below. One has $mpl(X,r_X)=2.$

Let $(Y, r_Y)$ be the one element solution, with $Y=\{ a \}.$ We
will build an extension $(Z,r) = X \natural Y.$ Since $Y$ is one
element solution, the left action on $Y$ on $X$ has to be via an
automorphism. Consider the permutation $\tau=(bc)(x_1x_2) \in
Sym(X).$ It is  easy to see that is an automorphism of $(X,r_X).$
Now one defines strong twisted union $(Z,r)=X\natural Y$ via the
action $\Lcal_{a}=(bc)(x_1x_2).$ We obtain  again the solution of
Example \ref{examplempl4}. Now one has $mpl(Z,r)=3 = max\{ mpl(X),
mpl(Y) \} + 1.$
\end{example}

\begin{lemma}
\label{graphlemma2} Let $x, a,b \in X,$ $\Lcal_a^x=(x_1 x_2\cdots
x_k)$, where $ x_1=x, k\geq 2,$ and let ${}^bx=x_{21} \neq x.$
Suppose that the left actions of $a$ and $b$ commute on $\Gamma_x$.
Then
\begin{enumerate}
\item
$\Lcal_a^{({}^bx)}=(x_{21} x_{22} \cdots x_{2k})$ is a cycle of
length exactly $k.$
\item If
 ${}^bx=x_{m+1}$ for some $m, 1 \leq m\leq k-1$, then
 $\Lcal_a^{{}^bx}=\Lcal_a^x,$ and $\Lcal_b^x =
(\Lcal_a^x)^m$.
 \item If part (2) does not apply then
 $\Lcal_a^x$  and $\Lcal_a^{{}^bx}$ are two disjoint
cycles.
\end{enumerate}
\end{lemma}
\begin{proof}
 We apply Lemma \ref{graphcellslemma} to
the commuting edges
\[
\xymatrix{
x_1 \ar[r]^a \ar[d]_ {b} & x_2  \\
x_{21}  & }
\]
and obtain the ``cell" subgraph   of {\bf $\Gamma$}:
\begin{equation}\label{com1}
\xymatrix{
x_1 \ar[r]^a \ar[d]_ {b} & x_2 \ar[d]^{b} \\
x_{21} \ar[r]_{a} & x_{22} .}
\end{equation}

Now we use this method to recursively build  the diagram
\begin{equation}\label{com2}
\xymatrix{ x_1 \ar[r]^{a} \ar[d]_{b} & x_2 \ar[r]^{a} \ar[d]_{b} &
\cdots \ar[r]^{a} & x_k \ar[r]^{a} \ar[d]_{b}&
x_1 \ar[r]^{a} \ar[d]_{b} & \cdots \\
 x_{21} \ar[r]^{a} \ar[d]_{b} & x_{22} \ar[r]^{a} \ar[d]_{b} &
\cdots \ar[r]^{a} & x_{2k} \ar[r]^{a} \ar[d]_{b}&
x_{21} \ar[r]^{a} \ar[d]_{b} & \cdots \\
\vdots &\vdots & \cdots  &  \vdots & \vdots &\cdots
}
\end{equation}
It is clear from the first two rows of the diagram  that the cycle
$\Lcal_a^{x_{21}}$ (represented by the second row) has length at
most $k$. Assume that its length $q$ is strictly less than $k$
then we obtain from (\ref{com2})
\begin{equation}\label{com3}
\xymatrix{ x_1 \ar[r]^{a} \ar[d]_{b} & x_2 \ar[r]^{a} \ar[d]_{b} &
\cdots \ar[r]^{a} & x_{q} \ar[r]^{a} \ar[d]_{b}&
x_{q+1} \ar[r]^{a} \ar[d]_{b} & \cdots \\
 x_{21} \ar[r]^{a} \ar[d]_{b} & x_{22} \ar[r]^{a} \ar[d]_{b} &
\cdots \ar[r]^{a} & x_{2q} \ar[r]^{a} \ar[d]_{b}& x_{21}
\ar[r]^{a}
\ar[d]_{b} & \cdots\\
\vdots &\vdots & \cdots  &  \vdots & \vdots &\cdots }
\end{equation}
which implies ${}^b{x_1}=x_{21}$ and ${}^b{x_{q+1}}=x_{21}.$ The
last is impossible in view of the nondegeneracy of $(X,r)$  and
$x_1 \neq x_{q+1}. $ We have shown that $\Lcal_a^{({}^bx)}$ is a
cycle of length exactly $k.$

Assume now ${}^bx=x_{21}=x_{m+1},$ for some $m\geq 1.$ Then we
deduce from  (\ref{com3}) the following diagram:
\[
\xymatrix{ x_1 \ar[r]^{a} \ar[d]_{b} & x_2 \ar[r]^{a} \ar[d]_{b} &
\cdots \ar[r]^{a} & x_{m+1} \ar[r]^{a} \ar[d]_{b}& \cdots
\ar[r]^{a} & x_k \ar[r]^{a} \ar[d]_{b}& x_1 \ar[r]^{a} \ar[d]_{b}
& \cdots
\\
 x_{m+1} \ar[r]^{a} \ar[d]_{b} & x_{m+2} \ar[r]^{a} \ar[d]_{b} &
\cdots \ar[r]^{a} & x_{2m+1} \ar[r]^{a} \ar[d]_{b}& \cdots
\ar[r]^{a} & x_{m+k} \ar[r]^{a} \ar[d]_{b}&
x_{m+1} \ar[r]^{a} \ar[d]_{b} & \cdots \\
x_{2m+1} \ar[r]^{a} \ar[d]_{b} & x_{2m+2} \ar[r]^{a} \ar[d]_{b} &
\cdots \ar[r]^{a} & x_{3m+1} \ar[r]^{a} \ar[d]_{b}& \cdots
\ar[r]^{a} & x_{2m+k} \ar[r]^{a} \ar[d]_{b}&
x_{2m+1} \ar[r]^{a} \ar[d]_{b} & \cdots \\
\vdots \ar[d]_{b}&\vdots \ar[d]_{b}&
 \cdots  &  \vdots \ar[d]_{b}
&
 \cdots  &  \vdots \ar[d]_{b}&
\vdots \ar[d]_{b}&\cdots\\
x_{sm+1} \ar[r]_ {a} \ar[d]_{b} & x_{sm+2} \ar[r]_ {a} \ar[d]_{b}
&\cdots \ar[r]_ {a} & \bullet \ar[r]_ {a} \ar[d]_{b}&\cdots
\ar[r]_ {a} & x_{sm+k} \ar[r]_ {a}\ar[d]_{b} & x_{sm+1}
\ar[r]_{a} \ar[d]_{b}& \cdots \\
\vdots &\vdots & \cdots  &  \vdots & \cdots & \vdots & \vdots
&\cdots }
\]
Looking at the  first column of the diagram we see that $\Lcal_b^x =
(\Lcal_a^x)^m.$ It follows from (\ref{remarkoncycles}) that
$\Lcal_a^{{}^bx}=\Lcal_a^{x_{m+1}} = \Lcal_a^{x}.$
This proves (1). The rest follows from (1).
\end{proof}
\begin{corollary}
\label{graphcor1} Let $x,  a, b \in X$ and suppose that the left
actions of $a$ and $b$ commute on $\Gamma_x$.
\begin{enumerate}
\item
If ${}^ax = {}^bx $ then $\Lcal_b^{x}=\Lcal_a^{x}.$
\item
If ${}^b{({}^ax)} = x$ then $\Lcal_b^{x}=(\Lcal_a^{x})^{-1}.$
\end{enumerate}
\end{corollary}
\begin{proof}
There is nothing to prove if ${}^ax =x.$ So suppose ${}^ax =
x_2\neq x.$ Let
 $\Lcal_a^{x}=(x_1 x_2\cdots x_k), k \geq 2,$ where, as usual, we set $x_1=x.$ Suppose ${}^ax =
 {}^bx.$ Then
\[
\xymatrix{
x_1 \ar[r]^a \ar[d]_ {b} & x_2  \\
x_2  & }
\]
is a subgraph of  $\Gamma$. The hypothesis of  Lemma
\ref{graphlemma2} part (2) is satisfied (here $m=1$). Therefore
$\Lcal_b^{x}=\Lcal_a^{x}.$ This proves part (1). Assume now
${}^b{({}^ax)} = x.$ This implies that
\[
\xymatrix{
x_1 \ar[r]^a  & x_2   \ar[d]_{b}  \\
  & x_1}
\]
is a subgraph of  $\Gamma$. It follows then  by Lemma
\ref{graphcellslemma} that $\Gamma$ contains also the following
subgraph:
\[
\xymatrix{
x_1 \ar[r]^a \ar[d]_ {b} & x_2 \ar[d]^{b} \\
x_{k} \ar[r]_{a} & x_{1}.}
\]
Apply Lemma \ref{graphlemma2} part (2) again to obtain
$\Lcal_b^{x}=(\Lcal_a^{x})^{k-1}=(\Lcal_a^{x})^{-1}$. (In the notation of the Lemma,
this time $m=k-1$).
\end{proof}
Clearly, in the case when $\Lcal_a^{x}=(x_1 x_2)$ the diagram in the
proof of part (2) of the corollary  becomes simply
\[
\xymatrix{
x_1 \ar[r]^a \ar[d]_ {b} & x_2 \ar[d]^{b} \\
x_{2} \ar[r]_{a} & x_{1}.}
\]
which in the real graph $\Gamma$ deforms to $x_1{\buildrel a,b \over
\longleftrightarrow} x_2.$

\begin{lemma}
\label{graphlemma3} Let $a,b , x \in X,$ let
$\Lcal_a^x=(x_1
x_2\cdots x_k)$, $\Lcal_b^x=(x_1 y_2\cdots y_s)$,
where $x=x_1$ and
$s \leq k.$ Suppose that the left actions of $a$ and
$b$ commute on
$\Gamma_x$ and suppose that the cycles $\Lcal_a^x$ and
$\Lcal_b^x$
have at least two common elements, i.e. $x_i=y_j$ for
some $i,j \neq
1$.
 Then
$\Lcal_b^x=(\Lcal_a^x)^m$, where $1 \leq m \leq k-1$.
\end{lemma}
\begin{proof} We claim that ${}^bx \in \{ x_2\cdots x_k \}.$
Assume the contrary, then we have $ \Lcal_b^x=(x_1y_{21}\cdots y_{p1}
x_{m+1}\cdots ),$ where $p \geq 1$ and $y_{21},\cdots ,y_{p1}$ are
distinct from $x_2,\cdots ,x_k.$

We want to deduce the vertices of the following subgraph of {\bf
$\Gamma$}:
\begin{equation}
\label{subgraph0} \xymatrix{ x_1 \ar[r]^{a} \ar[d]_{b} & x_2
\ar[r]^{a} \ar[d]_{b} & \cdots \ar[r]^{a} & x_k \ar[r]^{a}
\ar[d]_{b}&
x_1 \ar[r]^{a} \ar[d]_{b} & \cdots \\
  y_{21} \ar[r]^{a} \ar[d]_{b} & \bullet \ar[r]^{a} \ar[d]_{b} &
\cdots \ar[r]^{a} & \bullet \ar[r]^{a} \ar[d]_{b}&
\bullet \ar[r]^{a} \ar[d]_{b} & \cdots \\
\vdots \ar[d]_{b}&\vdots \ar[d]_{b}& \cdots  &  \vdots \ar[d]_{b}&
\vdots \ar[d]_{b}&\cdots\\
y_{p1} \ar[r]_ {a} \ar[d]_{b} & \bullet \ar[r]_ {a} \ar[d]_{b} &
\cdots \ar[r]_ {a} & \bullet \ar[r]_ {a}\ar[d]_{b} & \bullet
\ar[r]_{a} \ar[d]_{b}& \cdots \\
x_{m+1} \ar[r]_ {a} \ar[d]_{b} & \bullet \ar[r]_ {a} \ar[d]_{b} &
\cdots \ar[r]_ {a} & \bullet \ar[r]_ {a}\ar[d]_{b} & \bullet
\ar[r]_{a} \ar[d]_{b}& \cdots \\
\vdots &\vdots & \cdots  &  \vdots & \vdots &\cdots }
\end{equation}
So starting with $i=2,$ we apply Lemma
\ref{graphlemma2} to successively deduce the vertices of $i$-th row,
$i= 2, 3, \cdots.$ This way we obtain
\begin{equation}
\label{subgraph1} \xymatrix{ x_1 \ar[r]^{a} \ar[d]_{b} & x_2
\ar[r]^{a} \ar[d]_{b} & \cdots \ar[r]^{a} & x_k \ar[r]^{a}
\ar[d]_{b}&
x_1 \ar[r]^{a} \ar[d]_{b} & \cdots \\
 y_{21} \ar[r]^{a} \ar[d]_{b} & y_{22} \ar[r]^{a} \ar[d]_{b} &
\cdots \ar[r]^{a} & y_{2k} \ar[r]^{a} \ar[d]_{b}&
y_{21} \ar[r]^{a} \ar[d]_{b} & \cdots \\
\vdots \ar[d]_{b}&\vdots \ar[d]_{b}& \cdots  &  \vdots \ar[d]_{b}&
\vdots \ar[d]_{b}&\cdots\\
y_{p1} \ar[r]_ {a} \ar[d]_{b} & y_{p2} \ar[r]_ {a} \ar[d]_{b} &
\cdots \ar[r]_ {a} & y_{pk} \ar[r]_ {a}\ar[d]_{b} & y_{p1}
\ar[r]_{a} \ar[d]_{b}& \cdots \\
x_{m+1} \ar[r]_ {a} \ar[d]_{b} & x_{m+2} \ar[r]_ {a} \ar[d]_{b} &
\cdots \ar[r]_ {a} & x_{m+k} \ar[r]_ {a}\ar[d]_{b} & x_{m+1}
\ar[r]_{a} \ar[d]_{b}& \cdots \\
\vdots &\vdots & \cdots  &  \vdots & \vdots &\cdots }
\end{equation}
Note first that each row of (\ref{subgraph1}) represents a cycle of
length exactly $k$. Compare $\Lcal_a^x$ with $\Lcal_b^x$
(represented by the first column of (\ref{subgraph1})).
The hypothesis of Lemma~\ref{graphlemma2} is satisfied
but this time
we have that $\Lcal_a^x=(\Lcal_b^x)^{qp}$, where
$qp\equiv 2 (mod
k)$. Thus the lengths of the two cycles $k$ and $s$
satisfy  $k \leq
s.$ An equality is impossible, since by our assumption
the entry
$y_{21}$ of $\Lcal_b^x$ does not occur in $\Lcal_a^x.$
This yields
$k<s,$ a contradiction with the hypothesis. It follows
then that
${}^bx = x_{m+1}$ for some $m\geq 1$ as we claimed.
Now Lemma
\ref{graphlemma2} again implies
$\Lcal_b^x=(\Lcal_a^x)^m$.
\end{proof}

\begin{definition}
\label{Stardef} Let $(X,r)$ be a nondegenerate
involutive square-free solution of the YBE,
and $x \in X$. The set $Star(x)$ is defined as
\[
\Star(x)=\{ \Lcal_a^x \mid \; a \in X\}.
\]
\end{definition}

\begin{definition} By $\Gamma(x)$ we denote the
subgraph of (the
complete graph) $\Gamma(X,r)$ with a set of vertices
consisting of
all $x_j$ that occur in the cycles of $\Star(x),$ and
all edges
inherited from $\Gamma(X,r).$

On the set of stars we introduce an equivalence
relation $\approx$
defined as $\Star(x)\approx \Star(y)$ \emph{iff } the
complete
graphs $\Gamma(x)$ and $\Gamma(y)$ are isomorphic (as
labeled
graphs).
\end{definition}

As before, when we draw the graph $\Gamma(x)$ we shall
often present
it as a simplified graph in which if
$\Lcal_a^x=\Lcal_b^x,$ for $a
\neq b \in X$, we only draw it once. In various  cases
is convenient
to use even more schematic graph in which for each
$\Lcal_a^x=(x_1
x_2 \cdots x_k)\in \Star(x)$ with $x_1=x$ we draw
\[ \stackrel{x}{\bullet} \longrightarrow \bullet
\longrightarrow
\cdots \longrightarrow \bullet \] where each  variable
occurring in
$\Lcal_a^x$ is represented by a $\bullet$ (the last
corresponds to
$x_k$).

\begin{figure}
\begin{center}
\includegraphics[scale=.9]{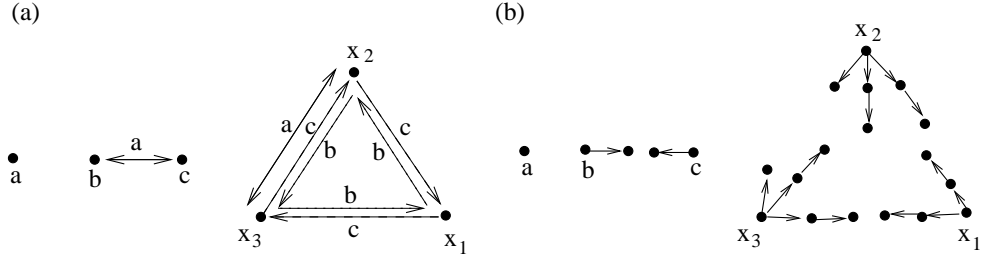}
\caption{Graph of Example~\ref{ejemploStar} illustrating $Star$ at each vertex.}
\label{ej5cap3}
\end{center}
\end{figure}

\begin{example}\label{ejemploStar}
Consider the nondegenerate square-free involutive solution $(X,r)$
defined on the
set
 \[
 X=\{x_1,x_2, x_3, a, b, c\},
 \]
 via the left actions
\[
\begin{array}{ll}
\Lcal_{x_1}=\Lcal_{x_2}=\Lcal_{x_3}=id_X,\quad &
\Lcal_{b}=(x_1x_2x_3),\\
\Lcal_{c}=(x_1x_3x_2),\ & \Lcal_{a}=(bc)(x_2x_3),
\end{array}
\]
The graph $\Gamma(X,r)$ is presented in figure
\ref{ej5cap3} (a). The graphs of the Stars of all elements of $X$ are
represented in part (b).
\end{example}

We introduce now  a particular type of subgraphs $\Gamma_1$ we call
them \emph{graphs of first type}.
\begin{definition}
\label{graph1} In our usual assumption and notation, let $\Gamma_1$
be a connected component of $\Gamma= \Gamma(X,r),  X_1 =
\Vcal(\Gamma_1)$ be its set of vertices.
 $\Gamma_1$ is \emph{a graph of first type}
if it contains a nontrivial edge and for each $a, b \in
\Ecal(\Gamma_1)$ the left actions of $a$ and $b$ commute on
$\Gamma_1$. In other words, if
\[
\Lcal_{({}^ba)\mid X_1} =\Lcal_{a\mid X_1}, \quad  \forall a,b
\in\Ecal(\Gamma_1).
\]
Next we define
\begin{equation}
\label{defG1} \Gcal_1=\Gcal(\Gamma_1): = {}_{gr} \langle
\Lcal_{a\mid X_1}; a \in \Ecal(\Gamma_1) \rangle.
\end{equation}
\end{definition}

\begin{lemma}
\label{graphlemma6} Let $\Gamma_1$ be a nontrivial connected
component of $\Gamma(X,r)$ with a set of labels $\Ecal_1=
\Ecal(\Gamma_1)$. Then $\Gamma_1$ is a graph of first type \emph{iff}
\begin{equation*}
 \Lcal_{(a^x)\mid \Gamma_1} =  \Lcal_{a\mid \Gamma_1}
\quad \forall a,x \in \Ecal_0.
\end{equation*}
\end{lemma}
\begin{proof}
The lemma follows from Definition \ref{graph1}, the cyclic
conditions and \textbf{lri}.
\end{proof}
Note that since $(X,r)$ is a solution, condition \textbf{l1} hods
and can be written as
\[
\Lcal_{a}\circ\Lcal_{b}= \Lcal_{{}^ab}\circ\Lcal_{a^b}\quad \forall
a,b \in X.
\]
It follows then from  Definition \ref{graph1} that the group
$\Gcal(\Gamma_1)$ is an abelian subgroup of $Sym(X_1).$ We recall
the following fact about finite abelian groups.
\begin{fact}
\label{fact1} \textbf{Basis Theorem}\cite[Ch. 2]{Sc} Let $G$ be a
finitely generated abelian group. If $s$ is the least integer such
that $G$ is generated by $s$ elements, then $G$ is the direct sum
(in abelian notation) of $s$ cyclic subgroups.
\end{fact}
\begin{remark}
\label{least_s} Note that, in general, the presentation of a finite
abelian group $G$ as a direct product of cyclic subgroups  is not an
invariant of the group, it becomes invariant if one considers only
direct products of cyclic subgroups of prime power order. (For
example one has $G= \mathbb{Z}_6=\mathbb{Z}_2 \times{\mathbb{Z}_3}$,
here $s(G)=1$ ).

However, the  least number $s=s(G)$ of a generating set of a finite
group is well determined, and we shall use it to introduce the
notion of \emph{basis of labels} for $\Gamma_1$. \end{remark}

It follows by Fact \ref{fact1}, that the group
$\Gcal_1=\Gcal(\Gamma_1)$ is the direct product of cyclic subgroups:
\begin{equation}
\label{basis1}\Gcal_1= (\Lcal_{a_1\mid X_1})\times (\Lcal_{a_2\mid
X_1})\times\cdots\times (\Lcal_{a_s\mid X_1}),
\end{equation}
where $s= s(\Gamma_1)$ is  the least integer such that $\Gcal_1$
is generated by $s$ elements. Note that $s$ is also an invariant of
$\Gamma_1$.
\begin{definition}
\label{basis2} A  set of labels $\mathbb{B}= \{a_1,\cdots ,a_s\}
\subseteq \Ecal(\Gamma_1)$ is \emph{a basis of of labels for
$\Gamma_1$} if (\ref{basis1}) holds and $s$ is the least number of
a generating set of  $\Gcal_1$.
\end{definition}

For computations it is convenient to use the arrows in $\Gamma_1$
labeled by $a_i, 1 \leq i \leq s$ as ``axes" of an $s$-dimensional
discrete "coordinate system". Let $k_i = k(a_i)$ as in Theorem
\ref{graphprop1} part (4).   We can chose an arbitrary vertex $x$ of
$\Gamma_1$ as an ``origin". It has coordinates $x=(0\cdots 0).$ Then
each vertex $y\in X_1$ has $s$ ``coordinates" which we define as
follows.
\[
\Lcal_{a_{1}}(00\cdots 0)= (10\cdots 0),\; \Lcal_{a_{2}}(00\cdots
0)= (01\cdots 0),\; \cdots, \; \Lcal_{a_{s}}(00\cdots 0)= (00\cdots
1).
\]
Clearly, every vertex $y$ of $\Gamma_1$ can be obtained by applying
a finite sequence of left actions:
\begin{equation}
\label{eqcoordinates1} y = (\Lcal_{a_{1}})^{m_1} \circ
(\Lcal_{a_{2}})^{m_2}\circ\cdots \circ
(\Lcal_{a_{s}})^{m_s}((0\cdots 0)) = (m_1 m_2  \cdots m_s),
\end{equation}
where $0 \leq m_i\leq k_i$.

We will now give a detailed description of the shape of a graph of
first type.
\begin{theorem}
\label{graphprop1} Let $(X,r)$ be a finite nondegenerate
square-free symmetric set and suppose the connected component
$\Gamma_1$ of $\Gamma(X,r)$ is a graph of first type. Let
$X_1=\Vcal(\Gamma_1),$ and $\Gcal_1=\Gcal(\Gamma_1)$, as in
(\ref{defG1}). Then the following conditions hold.
\begin{enumerate}
\item
\label{it1} If the nontrivial edge $x {\buildrel a \over
\longrightarrow}y, x \neq y$ occurs in $\Gamma_1$ then for every
vertex $z\ in \Vcal(\Gamma_1)$ there is an edge  $z  {\buildrel a
\over \longrightarrow} t, z \neq t\in \Vcal(\Gamma_1).$
\item
\label{it2} $X_1$ is $r$-invariant, and $(X_1, r_1)$ with $r_1=
r_{\mid X_1\times X_1}$ is a trivial solution.
\item
\label{it3}
Each $a \in X$ acts on $X_1$ as an automorphism, and
$\Gcal_1=\Gcal(\Gamma_1)$ is a normal subgroup of $Aut(X_1,r_1).$
\item
\label{it4} For every $a \in \Ecal(\Gamma_1),$ there exists  an
integer $k=k(a,\Gamma_1)$, which is  an invariant of $a$ and
$\Gamma_1$ such that for every vertex $x\in X_1 $ the cycle
$\Lcal_a^x$ is of length $k.$ Moreover
\[
\Lcal_{a\mid X_1 }= (x_1\cdots x_k) (y_1\cdots y_k)\cdots (z_1\cdots
z_k),
\]
where each vertex of $\Gamma_1$ occurs exactly once and the product
contains all disjoint cycles with edges labeled $a.$
\item
\label{it5} All vertices $x_1,\cdots ,x_N \in X_1$ have equivalent
stars, $Star(x_1) \approx Star(x_i), 1\leq i \leq N$.
\item
\label{it6} Let $\mathbb{B}= \{a_1,\cdots ,a_s\}$ be a basis of
labels for $\Gamma_1,$ and let $k_i = k(a_i,\Gamma_1), 1\leq i \leq
s,$ as in \ref{it4}. Then $\Gamma_1$ is a "multi-dimensional cube"
of dimension $s,$ and the order of $X_1$ (i.e. the size of the cube)
is exactly the product $N= k_1k_2\cdots k_s.$
\end{enumerate}
\end{theorem}
\begin{proof}
We start with part (\ref{it1}). By hypothesis $\Gamma_1$ is a
connected component of $\Gamma,$ and all its edges commute. Then it
follows from Lemma \ref{graphcellslemma} that for any ``neighbour"
$z$ of $x$, ($z=y$ is also possible) with
\[
\xymatrix{ x \ar[r]^a \ar[d]_ {b} & y\\ z & }
\]
 $\Gamma_1$ contains the subgraph:
\begin{equation}
\label{lemgraph1} \xymatrix{
x \ar[r]^a \ar[d]_ {b} & y \ar[d]^{b} \\
z \ar[r]_{a} & t .}
\end{equation}
Furthermore, by  Lemma \ref{graphcellslemma} part (3), $x \neq y$ implies
$z \neq t,$ which verifies (\ref{it1}).

\textbf{(\ref{it2})}. $X_1$ is the set of vertices of a connected
component of $\Gamma(X,r)$, so by Lemma \ref{graphlemma1} it is
$r$-invariant. Let $r_1= r_{\mid X_1\times X_1}$, we will show that
$(X, r_1)$ is a trivial solution, that is $r(x,y) =(y,x),$ for all
$x, y \in X_1.$ In view of \textbf{lri }  this is equivalent to
${}^yx=x,$ for all $x,y \in X_1$. Assume the contrary, there exist
$x,y \in X_1,$ such that ${}^yx \neq x.$ Then $\Gamma_1$ contains
the nontrivial  edge $x {\buildrel y \over \longrightarrow}\xi.$
This by (\ref{it1}) implies that every vertex $z \in X_1$ has an
edge $z {\buildrel y \over \longrightarrow}t, t \neq z.$ In
particular $y {\buildrel y \over \longrightarrow} \zeta, \zeta \neq
y,$ which is impossible, since $(X,r)$ is a square-free solution,
thus ${}^yx = y, \forall x,y \in X$. This implies (\ref{it2}).

The lemma below proves part (\ref{it3})  and will be needed also for
the proof of other statements in the sequel. It remains to prove the
last three items of the theorem. Suppose $x \in X_1, a,b \in
\Ecal(\Gamma_1),$ then by  Lemma \ref{graphlemma2} the cycles
$\Lcal_a^x$ and $\Lcal_a^{{}^bx}$ are of the same length, say $k$.
This and part (\ref{it1}) of the theorem  imply that for every
vertex $z\in X_1$ the cycle $\Lcal_a^z$ has length $k$ which proves
(\ref{it4}). Clearly, $(\ref{it1})\; \text{and}\;(\ref{it4})
\Longrightarrow (\ref{it5}) \Longrightarrow (\ref{it6}).$
\end{proof}

\begin{lemma}
\label{graphlemma5}  Under the hypothesis of Theorem
\ref{graphprop1},
\begin{enumerate}
\item
\label{item1}
The following equalities hold:
\[ {}^{a^y}x={}^ax ={}^{{}^ya}x  \quad \forall x,y\in
X_1, \; \text{and}\; \forall  a \in X.\]
\item
\label{graph1aut} Each $a \in X$ acts on $X_1$ as an automorphism,
and $\Gcal(\Gamma_1)$ is a normal subgroup of $Aut(X_1,r_1).$
\end{enumerate}
\end{lemma}
\begin{proof}
 Suppose $x,y \in X_1, a \in X.$  We use part (\ref{it2})
 of the theorem and
\textbf{l1} on $(X, r)$ to obtain
\[
{}^ax=^{(\ref{it2})}\;{}^a({}^yx)=^{\textbf{l1}}\;
{}^{{}^ay}{({}^{a^y}x)}=^{(\ref{it2})}\;{}^{a^y}x.
\]
This imply ${}^{a^y}x={}^ax,$ for every $x,y\in X_1,$ and $a \in X,$
which verifies the left hand side of the equality of part~(\ref{item1}).
From there, using \textbf{lri} and the non-degeneracy of $r$ one
easily deduces the righthand side of the required equality.

It follows from part (\ref{item1}) and Lemma \ref{criterionLa} that
\[
\Lcal_{a\mid \Gamma_{1}}\in Aut(X_1,r_1), \forall a \in X
\]
thus $\Gcal(\Gamma_1)$ is a subgroup of $Aut(X_1,r_1).$ Now Lemma
\ref{leftright}, and the equality (\ref{autre}) imply that
$\Gcal(\Gamma_1)$ is a normal subgroup of $Aut(X_1,r_1).$ The lemma
has been proved, which also completes the proof of the theorem.
\end{proof}

\begin{theorem}
\label{mpl2theorem} Let $(X,r)$ be a finite nondegenerate
square-free symmetric set of order $\geq 2$, $G=G(X,r), \Gcal=
\Gcal(X,r),$ $\Gamma=\Gamma(X,r)$, Aut(X,r) in the usual notation.
Let $\Gamma_1, \Gamma_2, \cdots, \Gamma_s$ be all connected
components, and $X_i=\Vcal(\Gamma_i), 1 \leq i \leq s,$
respectively, be their sets of vertices. The following conditions
are equivalent.
\begin{enumerate}
\item
\label{mpl2it1}
 $(X,r)$ is multipermutation solution with $mpl(X,r)=2.$
\item
\label{mpl2it2} $\Gcal(X,r)$ is an abelian group of order $\geq 2$.
\item
\label{mpl2it3} $\{id_X\}\neq \Gcal(X,r)$ is a subgroup of the
automorphism group $Aut(X,r)$.
\item
\label{mpl2it4} The set of nontrivial components is nonempty.
Suppose these are $\Gamma_1, \Gamma_2, \cdots, \Gamma_p,$ with $ 1
\leq p \leq s.$ Every nontrivial connected component $\Gamma_i$ is a
graph of first type. Furthermore, in this case, for each pair $i, j,
 1 \leq i,j \leq s$ and each $a,b \in X_j$  one has $\Lcal_{a\mid
X_i} = \Lcal_{b\mid X_i}.$
\item
\label{mpl2it5} $(X,r)$ can be split into disjoint $r$-invariant
subsets $X_i, 1 \leq i \leq s,$ where each $(X_i, r_{\mid X_i})$ is
a trivial solution and $X= X_1\natural X_2\natural \cdots \natural
X_s,$ in the sense that we can put parentheses and "apply"
$\natural$ in any order. In particular, for any pair $i,j$ $r$
induces $X_i\natural X_j$  and
\begin{equation}
\label{mpl2eq0} X = X_i\natural (\bigcup_{1 \leq j \leq s, j\neq
i}X_j)
\end{equation}
\end{enumerate}
\end{theorem}
\begin{proof}
We shall show that \[ (\ref{mpl2it1})\Longleftrightarrow
(\ref{mpl2it2})\Longleftrightarrow (\ref{mpl2it3})\; \text{and} \;
(\ref{mpl2it2})\Longleftrightarrow (\ref{mpl2it4})\Longrightarrow
(\ref{mpl2it5})\Longrightarrow(\ref{mpl2it2}).
\]
 Under the hypothesis of the theorem we start with
the following easy lemma.
\begin{lemma}
\label{mpl2lemma1} The following conditions are equivalent.

(i)  $mpl(X,r)=2$

(ii) $mpl(Ret(X,r)) =1,$

(iii) $Ret(X,r) = ([X], r_{[X]})$ is the trivial solution.

(iv) $ {}^{[x]}{[a]} = [{}^xa]= [a]$, for all $a,x \in X$

(v) There is an equality:
\begin{equation}
\label{mpl2eq1}
 \Lcal_{{}^xa} = \Lcal_{a} \quad \forall x, a \in X
\end{equation}
(vi) There are equalities
\begin{equation}
 \label{mpl2eq1a} \Lcal_{{}^xa}= \Lcal_{a}=
\Lcal_{a^x},  \quad \forall x,a \in X.
\end{equation}
\end{lemma}
\begin{proof}
It follows from Definition \ref{mpldef} and Lemma
\ref{mplretractlemma1} that (i) $\Longleftrightarrow$ (ii). Lemma
\ref{mpl1lem} gives the equivalence (ii) $\Longleftrightarrow$
(iii).  For (iii) $\Longleftrightarrow$ (vi), the condition that the
retract $([X], r_{[X]})$ is the trivial solution, implies that each
pair $x,a \in X$ satisfies
\[
{}^{[x]}{[a]} = [{}^xa]= [a]=[a]^{[x]}=[a^x]
\]
These equalities, ``translated" in terms of the left actions give
(\ref{mpl2eq1a}). (iii) $\Longleftrightarrow$ (iv) follows from
Definition \ref{actionson[X]def} and Lemma \ref{mpl1lem} again. By
Definition \ref{[x]def}
 the equality (\ref{mpl2eq1}) is just
 expression of
$[{}^xa]= [a]$ in terms of the left action, thus (iv)
$\Longleftrightarrow$ (v).
\end{proof}
We continue with the proof of the theorem.
 $(\ref{mpl2it1})\Longrightarrow(\ref{mpl2it2})$. Suppose $mpl X = 2$.
 Recall that
$\Gcal(X,r)$ is the subgroup of $Sym(X)$ generated by the
permutations $\Lcal_{a}, a\in X$. By Lemma \ref{mpl2lemma1} the
equalities (\ref{mpl2eq1a}) are in force, so we yield:
\[
\Lcal_{a}\circ \Lcal_{b}=^{\textbf{l1}}\;\Lcal_{{}^ab}\circ
\Lcal_{a^b}=^{(\ref{mpl2eq1a})}\;\Lcal_{b}\circ \Lcal_{a},
\quad\forall a,b \in X
\]
Thus $\Gcal(X,r)$ is an abelian group.
The following implications
are clear
\[
 (\ref{mpl2it2})\Longrightarrow (\ref{mpl2eq1})
 \Longrightarrow^{\text{by Lemma }\ref{mpl2lemma1}} \quad (\ref{mpl2it1})
\]
The implications $(\ref{mpl2eq1})\; \Longleftrightarrow\;
(\ref{mpl2it3})$ follow from Lemma  \ref{criterionLa}, and together
with  Lemma \ref{mpl2lemma1}  yield $(\ref{mpl2it1})\;
\Longleftrightarrow\; (\ref{mpl2it3}).$

$(\ref{mpl2it2})\Longleftrightarrow (\ref{mpl2it4})$. Suppose
(\ref{mpl2it2}) holds. Then $\Gcal(X,r)\neq id_X,$  so there is a
nontrivial connected component of $\Gamma.$ In addition by
(\ref{mpl2eq1}) each such a component $\Gamma_i$ is a graph of first
type.

\begin{lemma}
\label{mpl2lemma2} $(\ref{mpl2it4})\; \Longrightarrow
\;(\ref{mpl2eq1})$.
\end{lemma}
\begin{proof}
Suppose  (\ref{mpl2it4}) holds and assume the contrary, there exist
a pair $a, x \in X,$ such that $\Lcal_{{}^xa} \neq \Lcal_{a}.$

Hence there exist a $t,$ for which
\begin{equation}
\label{mpl2eq2} {}^{{}^xa}t \neq {}^at.
\end{equation}
 Clearly, then the connected component
$\Gamma_t$ which contains $t$ is nontrivial. Denote by $X_0$ the set
of vertices of $\Gamma_t,$ the set of labels we denote by  $\Ecal_0
= \Ecal(\Gamma_t).$ By our assumption $\Gamma_t$ is a graph of first
type, and therefore all its edges commute, so
at least one of  $a$ and $x$ is not in $\Ecal_0$. Suppose $a$ is
not in  $\Ecal_0.$ Then ${}^{{}^xa}t\neq {}^at =t,$ implies that
\begin{equation}
\label{mpl2eq3} {}^xa \in \Ecal_0.
\end{equation}
Two cases are possible. \textbf{a)} $x \in \Ecal_0.$ Then, since
$\Gamma_t$ is a graph of first type, $\Lcal_{{}^xa}$ and $\Lcal_{x}$
commute on $\Gamma_t.$ Now the equalities
\begin{equation}
\label{mpl2eq4}
 \Lcal_{({}^xa\mid \Gamma_1)} =^{\text{ Lemma
}\ref{graphlemma6}} \;\Lcal_{(({}^xa)^x\mid \Gamma_1)}
=^{\text{\textbf{ lri }} }\; \Lcal_{a \mid \Gamma_1}
\end{equation}
give a contradiction with (\ref{mpl2eq2}). \textbf{b)} $x$ is not in
$\Ecal_0,$  so ${}^xt= t,$ and ${}^x{({}^at)}= {}^at.$ By
\textbf{l1} one has ${}^x{({}^at)}= {}^{({}^xa)}{({}^{(x^a)}t)},$ so
(\ref{mpl2eq2}) imply ${}^{(x^a)}t\neq t.$ Thus $x^a\in \Ecal_0,$
and by  (\ref{mpl2eq3}) Lemma \ref{graphlemma6} one has
\begin{equation}
\label{mpl2eq5} \Lcal_{({}^xa \mid \Gamma_1)}  =^{\text{by Lemma
}\ref{graphlemma6}} \Lcal_{(({}^xa)^{x^a}\mid \Gamma_1)} =
\; \Lcal_{a \mid \Gamma_1}.
\end{equation}
For the right hand side of (\ref{mpl2eq5}) we used the following
equalities which come from the cyclic condition {\rm\bf cl1}, see
Definition \ref{cyclicconditionsall}, and \textbf{lri}:
\[
({}^xa)^{x^a}=^{\text{by  \textbf{cl1}}}\;=
({}^{x^a}a)^{x^a}=^{\text{by \textbf{lri}}} \; a.
\]
This way (\ref{mpl2eq5}) gives a contradiction with (\ref{mpl2eq2}).
It follows then that $a \in\Ecal_0.$ Next (\ref{mpl2eq2}) implies
that $x$ is not  in $\Ecal_0,$ so it commutes with every element of
$X_0$ and in particular with ${}^at.$ Hence one has
\[
{}^at= {}^x{({}^at)}={}^{({}^xa)}{({}^{(x^a)}t)},
\]
and therefore $x^a \in \Ecal_0.$ It follows then that
\[
 \Lcal_{a \mid
\Gamma_1} = \Lcal_{(({}^{x^a}a)\mid \Gamma_1)}=^{\text{by
\textbf{cl1}}}\; \Lcal_{({}^xa \mid \Gamma_1)}.
\]
(Remind that due to \textbf{cl1 } one has ${}^{x^a}a= {}^xa,$ which
was used in the right hand side of the above equality). It follows
than that (\ref{mpl2eq2}) is impossible. We have shown that for each
$a,x \in X$ one has $\Lcal_{{}^xa}= \Lcal_{a}$. This proves the
lemma
\end{proof}

It follows from Lemmas \ref{mpl2lemma1} and \ref{mpl2lemma2} that
$(\ref{mpl2it4})\; \Longrightarrow \;(\ref{mpl2it2})$

Next we show $(\ref{mpl2it4})\Longrightarrow (\ref{mpl2it5})$. So,
assume now (\ref{mpl2it4}) and consider the sets $X_i,1 \leq
 i \leq s.$ By Theorem \ref{graphprop1}  $(X_i, r_i)$ is a trivial solution of order
 $\geq 2$ for all
 $i, 1 \leq
 i \leq p$ and (in case that $p < s$) it is one element solution for $i, p+1 \leq
 i \leq s$.
 by Lemma \ref{mpl2lemma2} the equalities (\ref{mpl2eq1}) hold, and it is
 easy to see that (\ref{mpl2eq1}) implies (\ref{mpl2it5}).

Conversely, suppose (\ref{mpl2it5}) holds. Let $a \in X.$ Then there
exist unique $i, 1\leq i \leq s$ such that $a\in X_i. $ It follows
then from (\ref{mpl2eq0}) that $\Lcal_{a}$ acts on $X \backslash X_i
$ as an automorphism.  Note that since $X_i$ is the trivial
solution, or one element solution thus in both cases one has
$\Lcal_{a\mid X_i}= id_{X_i},$ which clearly implies $\Lcal_{a} =
\Lcal_{a\mid X \backslash X_i}.$ It follows then that $\Lcal_{a} \in
Aut(X,r).$ We have shown $(\ref{mpl2it5})\Longrightarrow
(\ref{mpl2it2})$. The theorem has been proved.
\end{proof}

\begin{proposition}
\label{mpl3prop} Let $(X,r)$ be a nondegenerate
square-free
symmetric set. Let $\Gamma_1, \Gamma_2, \cdots,
\Gamma_s$ be all
connected components, and $X_i=\Vcal(\Gamma_i), 1 \leq
i \leq s,$
respectively, be their sets of vertices.  Suppose
$mpl(X, r) = 3.$
Then for each $i, 1 \leq i \leq s,$, $(X_i, r_i)$ with
$r_i =
r_{\mid X_i\times X_i}$ is a  multipermutation
solutions of level
$mpl(X_i)\leq 2.$  Firthermore, for each pair $i\neq
j, 1 \leq i,j
\leq s,$ one has $X_i\natural X_j.$
\end{proposition}
\begin{proof}
Proposition \ref{graphlemma1} implies that for each
$i, 1 \leq i
\leq s$ and each $x\in X_i$ one has $X_i \subseteq \Ocal(x, 2),$ and
therefore
$mpl(X_i,r_i) \leq 2.$

Chose now an arbitrary pair  $i\neq j, 1 \leq
i,j \leq s.$
We claim that $X_i\natural X_j.$ It will be enough to
show
\begin{equation}
\label{mpl3eq1} {}^{\alpha^y}x={}^{\alpha}x\quad
{}^{x^{\beta}}{\alpha}={}^{x}{\alpha}\quad\forall x,y
\in X_i,
\alpha, \beta\in X_j. \end{equation} Let $x,y \in X_i,$
$\alpha \in
X_j.$ Then ${\alpha}^{y} \sim {\alpha}^{x},$ so
\[
{}^{\alpha^y}x={}^{\alpha^x}x={}^{\textbf{cl1}}
\;{}^{\alpha}x,
\]
which proves the left hand side equality of
(\ref{mpl3eq1}). Similar
argument verifies the remaining equality in
(\ref{mpl3eq1}).
\end{proof}

We will illustrate our theory on Examples~\ref{autex1}
and~\ref{autex2}. To begin with,  it  follows  straightforwardly from the  definition of
$(X, r)$ in these examples that
for $Ret(X,r)=([X], r_{[X]})$ one has $[X]=
\{[x_1],[b],[c]\},$ and
$r_{[X])}$  the trivial solution. Therefore
$mpl(X)=2.$ Moreover,  $Z = X\natural\{a\},$ where $\Lcal_a=
\Lcal_b\circ\Lcal_c \in
Aut(X,r).$ We leave the reader to verify $mpl(Z,
r_Z)=3.$

Clearly the graph  $\Gamma(Z,r_Z)$ is obtained by adding to
$\Gamma(X,r)$ a new  vertex and  two new edges labeled by $a.$

Next we shall describe the group $Aut(Z,r_Z)$. Each
automorphism
$\varphi \in Aut(Z,r_Z)$ is a product
$\varphi=\varphi^{\prime}
\varphi^{\prime\prime}$
where $\varphi^{\prime} \in Sym(a,b,c),
\varphi^{\prime\prime}\in Sym(x_1, x_2, x_3, x_4) $. We claim that $\varphi$ is uniquely determined by the
data
$(\varphi^{\prime}, \varphi^{\prime\prime}(x_1)).$
Indeed, knowing the image
$\varphi^{\prime\prime}(x_1),$ one applies the
equalities
\[
\varphi\circ\Lcal_a=\Lcal_{\varphi(a)}\circ\varphi,\quad
\varphi\circ\Lcal_b=\Lcal_{\varphi(b)}\circ\varphi,\quad\varphi\circ\Lcal_c=
\Lcal_{\varphi(c)}\circ\varphi.
\]
to find $\varphi(x_i), 2\leq i \leq4.$ We can then tabulate the
automorphisms with row corresponding to $\varphi'\in S_3$ and column
corresponding to the value of $\varphi''(x_1)$, which in all except
the last row (where we use a different notation) this value supplies
the index used. For simplicity we write $1,2, \cdots $ instead of
$x_1, x_2,\cdots$. Then the full list of automorphisms is:
\[\begin{array}{llll} \tau_1=(bc)(23),& \tau_2=(bc)(1243),& \tau_3=(bc)(1342),& \tau_4=(bc)(14)\\
 \pi_1=(abc)(234),& \pi_2=(abc)(124),&\pi_3=(abc)(132),& \pi_4=(abc)(143)\\
 \eta_1=(acb)(243),& \eta_2=(acb)(123),& \eta_3=(acb)(134),& \eta_4=(acb)(142)\\
 \rho_1=(ab)(24),& \rho_2=(bc)(1234),&\rho_3=(ab)(13),& \rho_4=(ab)(1432)\\
\sigma_1=(ac)(34),& \sigma_2=(ac)(12),& \sigma_3=(ac)(1324),& \sigma_4=(ac)(1423)\\
id_Z=(),& \Lcal_b=(12)(34),& \Lcal_c=(13)(24) & \Lcal_a=(14)(23).\end{array}\]

 For example,  $\pi_2$
corresponds to the data $(\pi_{2 \mid\{a b c\}}=
(abc),
\pi_{2}(x_1)=x_2).$ We will do the first step of the
``search" for
the missing information to see how this was obtained.
\[
\pi_{2}\circ\Lcal_b(x_1)=\Lcal_{\pi_{2}(b)}\circ\pi_{2}(x_1),
\quad
\pi_{2}(x_2)=\Lcal_{c}(x_2)=x_4
\]
so $\pi_{2}$ acts as $x_1\mapsto x_2 \mapsto
x_4.$ Similar
computation with $\pi_{2}\circ\Lcal_a(x_4)$ shows that
$\pi_{2}(x_4) = x_1.$ Hence
$\pi_{2}=(abc)(x_1x_2x_4).$  Note that the $\eta$ row consists of the inverses of the $\pi$ row.  In the same way one obtains the rest of the table.
 This list was illustrated  in Figure~\ref{figautex2}.

We claim that $Aut(Z,r_z)$ is isomorphic to the
symmetric group
$S_4.$ Indeed, we know that $S_4$ is generated by
$\{(1234),
(12)\}.$ Direct computation shows that the assignment
\[
(1234)\longrightarrow (ab)(1234) =\rho_2,\quad (12)
\longrightarrow
(ac)(12)=\sigma_2
\]
extends to an isomorphism of groups $S_4
\longrightarrow Aut(X,r).$

The computation of $Aut(X,r)$ as a set is done
analogously, and the
relations show directly that there is a group
isomorphism
 $Aut(X,r) \approx D_4.$ Clearly,
$Aut(X,r)$ is a subgroup of
$ Aut(Z,r_Z)$ as expected.

\bibliographystyle{ams-alpha}

\end{document}